\documentclass[a4paper,12pt,reqno]{amsart}

\usepackage{amsmath,amsfonts,amsthm,amssymb,mathrsfs}
\usepackage{a4wide,amsmath,latexsym,euscript,units,color}
\usepackage[latin1]{inputenc}

\newcommand{\der}{\delta}

\newcommand{\hsi}{\hat \sigma}

\newcommand{\cal}{\mathcal}
\newcommand{\id}{\mbox{Id}}

\newcommand{\iot}{\int_{0}^{t}}
\newcommand{\iott}{\int_{0}^{T}}
\newcommand{\ist}{\int_{s}^{t}}
\newcommand{\norm}[1]{\lVert #1\rVert}

\newcommand{\ott}{[0,T]}

\newcommand{\xd}{{\bf x^{2}}}
\newcommand{\bd}{{\bf B^{2}}}

\newcommand{\bdst}{{\bf B}_{st}^{\bf 2}}

\newcommand{\1}{{\bf 1}}
\def\lpa{\langle}
\def\rpa{\rangle}

\newcommand{\R}{\mathbb R}
\newcommand{\N}{\mathbb N}

\newcommand{\cb}{\mathcal B}
\newcommand{\cac}{\mathcal C}

\newcommand{\ce}{\mathcal E}

\newcommand{\ch}{\mathcal H}

\newcommand{\cj}{\mathcal J}
\newcommand{\cl}{\mathcal L}

\newcommand{\cn}{\mathcal N}
\newcommand{\cq}{\mathcal Q}

\newcommand{\cw}{\mathcal W}
\newcommand{\cz}{\mathcal Z}

\newcommand{\al}{\alpha}
\newcommand{\ep}{\varepsilon}
\newcommand{\e}{\varepsilon}

\newcommand{\ga}{\gamma}
\newcommand{\ka}{\kappa}

\newcommand{\laa}{\Lambda}

\newcommand{\si}{\sigma}

\newcommand{\vp}{\varphi}

\newcommand{\lp}{\left(}
\newcommand{\rp}{\right)}
\newcommand{\lc}{\left[}
\newcommand{\rc}{\right]}
\newcommand{\lcl}{\left\{}
\newcommand{\rcl}{\right\}}

\newcommand{\lla}{\left\langle}
\newcommand{\rra}{\right\rangle}

\newtheorem{theorem}{Theorem}[section]

\newtheorem{corollary}[theorem]{Corollary}

\newtheorem{definition}[theorem]{Definition}

\newtheorem{hypothesis}[theorem]{Hypothesis}
\newtheorem{lemma}[theorem]{Lemma}

\newtheorem{proposition}[theorem]{Proposition}

\theoremstyle{remark}
\newtheorem{remark}[theorem]{Remark}

\begin{document}

\title{Weak approximation of a fractional SDE}
\author[X. Bardina \and I. Nourdin \and C. Rovira \and S. Tindel]{X. Bardina \and I. Nourdin \and C. Rovira  \and  S. Tindel}
\date{\today}

\begin{abstract}
In this note, a diffusion approximation result is shown for stochastic
differential equations driven by a (Liouville) fractional Brownian motion $B$ with Hurst
parameter $H\in(1/3,1/2)$. More precisely, we resort to the Kac-Stroock type approximation
using a Poisson process studied in \cite{BJT,DJ}, 
and our method of proof relies on the algebraic
integration theory introduced by Gubinelli in \cite{Gu}.
\end{abstract}



\address{
{\it Xavier Bardina:}
{\rm Departament de Matem\`atiques, Facultat de Ci\`encies, Edifici C, Universitat Aut\`onoma de Barcelona, 08193 Bellaterra, Spain}.
{\it Email: }{\tt Xavier.Bardina@uab.cat}
\newline
$\mbox{ }$\hspace{0.1cm}
{\it Ivan Nourdin:}
{\rm Laboratoire de Probabilit\'es et Mod\`eles Al\'eatoires,
Universit{\'e} Pierre et Marie Curie,
Bo{\^\i}te courrier 188, 4 Place Jussieu, 75252 Paris Cedex 5, France}.
{\it Email: }{\tt ivan.nourdin@upmc.fr}
\newline
$\mbox{ }$\hspace{0.1cm}
{\it Carles Rovira:}
{\rm Facultat de Matem\`atiques, Universitat de Barcelona, Gran Via 585, 08007 Barcelona, Spain}.
{\it Email: }{\tt carles.rovira@ub.edu}
\newline
$\mbox{ }$\hspace{0.1cm}
{\it Samy Tindel:}
{\rm Institut {\'E}lie Cartan Nancy, B.P. 239,
54506 Vand{\oe}uvre-l{\`e}s-Nancy Cedex, France}.
{\it Email: }{\tt tindel@iecn.u-nancy.fr}
}

\maketitle

\section{Introduction}

After a decade of efforts \cite{ALN,CQ,Gu,lejay,LyonsBook,nourdin-simon,nual-cours},
it can arguably be said that the basis of the
stochastic integration theory with respect to a rough path in general, and with
respect to a fractional Brownian motion (fBm) in particular,
has been now settled in a rather simple and secure way. This allows in particular to
define rigorously and solve equations on an arbitrary interval $\ott$ with $T>0$, of the form:
\begin{equation}\label{eq:eds-intro}
dy_t=\si\lp y_{t} \rp dB_t
+ b\lp   y_{t} \rp dt,\quad y_0=a\in\R^n,
\end{equation}
where $\sigma:\R^n\rightarrow \R^{n\times d}$, $b:\R^n\rightarrow\R^n$ are two
bounded  and smooth functions, and $B$ stands for a $d$-dimensional
fBm with Hurst parameter $H>1/4$. A question which
arises naturally in this context is then to try to establish some of the basic
properties of the process $y$ defined by (\ref{eq:eds-intro}), and this global program
has already been started as far as moments estimates \cite{HN}, large deviations \cite{LQZ,millet},
or properties of the law \cite{BC,NNRT} are concerned (let us mention at this point that the forthcoming book \cite{FVbk} will give a detailed account on most of these topics).

\vspace{0.2cm}

In the current note, we wish to address another natural problem related to the fractional
diffusion process $y$ defined by (\ref{eq:eds-intro}). Indeed, in the case where $B$ is an
ordinary Brownian motion, one of the most popular method in order to simulate $y$ is
the following:
approximate $B$ by a sequence of smooth or piecewise linear functions, say $(X^{\e })_{\e>0}$,
which converges in law to $B$,
e.g. an interpolated and rescaled random walk. Then see if the process $y^{\e}$ solution
of equation (\ref{eq:eds-intro}) driven by $X^{\e}$ converges in law, as a process, to $y$.
This kind of result, usually known as diffusion approximation, has been thoroughly
studied in the literature (see e.g. \cite{JLK,TZ,WZ}),
since it also shows that equations like (\ref{eq:eds-intro}) may
emerge as the limit of a noisy equation driven by a fast oscillating function.
The diffusion approximation program has also been taken up in the fBm case by Marty in
\cite{marty}, with some random wave problems in mind,
but only in the cases where $H>1/2$ or the dimension $d$ of the fBm is 1. 
Also note that, in a more general context, strong and weak approximations
to Gaussian rough paths have been studied systematically by Friz and Victoir in \cite{FV}. Among other results, 
the following is proved in this latter reference: let $(X^\ep)_{\ep >0}$ be a sequence of $d$-dimensional centered Gaussian 
processes with independent components and covariance function $R^\ep$. Let $X$ be another $d$-dimensional centered Gaussian 
processes with independent components and covariance function $R$. Assume that all those processes admit a rough path of 
order 2, that $R^\ep$ converges pointwise to $R$, and that $R^\ep$ is suitably dominated in $p$-variation norm for some 
$p\in[1,2)$. Then the rough path associated to $X^\ep$ also converges weakly, in $2p$-variation norm, 
to the rough path associated to $X$. 

\smallskip

This result does not close the diffusion approximation problem for solutions of SDEs like ~(\ref{eq:eds-intro}). 
Indeed, for computational and implementation reasons, the most typical processes taken as approximations to $B$ are non Gaussian, and more specifically, are usually 
based on random walks \cite{KP,WZ,So} or Kac-Stroock's type \cite{BJT,DJ,kac,stroock} approximations. However, the issue of diffusion approximations in a non-Gaussian context has hardly been addressed in the literature, and we are only aware of the aforementioned reference \cite{marty}, as well as the recent preprint \cite{FF} (which deals with Donsker's theorem in the rough path topology) for significant results on the topic. The current article proposes then a natural step in this direction, 
and studies diffusion approximations to (\ref{eq:eds-intro}) based on Kac-Stroock's approximation to white noise.

\smallskip

Let us be more specific about the kind of result we will obtain. First of all, we consider in the sequel the so-called $d$-dimensional Liouville fBm $B$, with Hurst parameter $H\in(1/3,1/2)$, as the driving process of equation (\ref{eq:eds-intro}). This is convenient for computational reasons (especially for the bounds we use on integration kernels), and is harmless in terms of generality, since the difference between the usual fBm and Liouville's one is a finite variation process (as shown in \cite{AMN}).
More precisely, we assume that $B$ can be written as $B=(B^1,\ldots,B^d)$, where the $B^i$'s are
$d$ independent centered Gaussian processes of the form
$$
B_t^i=\iot (t-r)^{H-\frac12} dW_r^i,
$$
for a $d$-dimensional Wiener process $W=(W^1,\ldots,W^d)$.
As an approximating sequence of $B$, we shall choose $(X^{\e})_{\e>0}$, where
$X^{\e,i}$ is defined 
as follows, for $i=1,\ldots,d$:
\begin{equation}\label{eq:def-x-ep}
X^{i,\e}(t)=\int_0^t (t+\e-r)^{H-\frac12}\theta^{\e,i}(r)dr,
\end{equation}
where
\begin{equation}\label{stroock}
\theta^{\e,i}(r)=\frac1{\e}(-1)^{N^i(\frac{r}{\e})},
\end{equation}
for $N^i$, $i=1,\ldots, d$, some independent standard Poisson processes.
Let us then consider the process $y^\e$ solution to equation (\ref{eq:eds-intro})
driven by $X^\e$, namely:
\begin{equation}\label{eq:eds-approx}
dy_t^\e=\si\lp y_{t}^\e \rp dX_t^\e
+ b\lp   y_{t}^\e \rp dt,\quad y_0^\e=a\in\R^n, \quad t\in\ott.
\end{equation}
Then our main result is as follows:
\begin{theorem}\label{thm:approx-diffusion}
Let $(y^\e)_{\e>0}$ be the family of processes defined by (\ref{eq:eds-approx}),
and let $1/3<\gamma<H$, where $H$ is the Hurst parameter of $B$.
Then, as $\ep\to 0$, $y^\ep$ converges in law to the process $y$ obtained as
the solution to (\ref{eq:eds-intro}), where the convergence takes place in the
Hölder space $\cac^\ga(\ott;\R^n)$.
\end{theorem}
Observe that we have only considered the case $H>1/3$ in the last result. This is of course for computational and notational sake, but it should also be mentioned that some of our kernel estimates, needed for the convergence in law, heavily rely on the assumption $H>1/3$. 
On the other hand, the case $H>1/2$ follows easily from the results contained in \cite{DJ}, 
and the case $H=1/2$ is precisely Stroock's result \cite{stroock}. This is why our future computations focus on the case $1/3<H<1/2$.

\smallskip

The general strategy we shall follow in order to get our main result is rather natural in
the rough path context: it is a well-known fact that the solution $y$ to (\ref{eq:eds-intro}) is a continuous function of
$B$ and of the L\'evy area of $B$ (which will be called $\bd$), considered as
elements of some suitable Hölder (or $p$-variation) spaces. Hence, in order to obtain the convergence $y^\ep{\to}y$ in law, it will be sufficient to check the
convergence of the corresponding approximations $X^\e$ and ${\bf X}^{\bf 2,\e}$ in their respective Hölder spaces (observe however that ${\bf X}^{\bf 2,\e}$ is not needed, in principle, for the definition of $y^\ep$). Then the two main technical problems we will have to
solve are the following:
\begin{enumerate}
\item
First of all, we shall use the {\it simplified} version of the rough path formalism,
called algebraic integration,
introduced by Gubinelli in \cite{Gu},
which will be summarized in the next section. In the particular context of
weak approximations, this allows us to deal with approximations of $B$ and $\bd$ directly,
without recurring to discretized paths as in \cite{CQ}. However, the
algebraic integration formalism relies on some space $\cac_k^\ga$, where
$k$ stands for a number of variables in $\ott$, and $\ga$ for a Hölder type exponent.
Thus, an important step will be to find a suitable tightness criterion in these spaces.
For this point, we refer to Section \ref{sec:tightness}.
\item The convergence of finite dimensional distributions (``fdd'' in the sequel) for the
L\'evy area $\bd$ will be proved in Section \ref{sec:fdd-cvgce}, and will be based on some sharp 
estimates concerning the Kac-Stroock kernel (\ref{stroock}) that are performed in Section \ref{sec:tech}. 
Indeed, this latter section is mostly devoted to quantify the distance between 
$\int_0^T f(u)\theta^\e(u)du$ and $\int_0^T f(u)dW_u$ for a smooth enough function $f$, 
in the sense of characteristic functions. This constitutes a generalization of \cite{DJ}, 
which is interesting in its own right.
\end{enumerate}

\vspace{0.2cm}

Here is how our paper is structured: in Section \ref{sec:basic-alg-intg}, we shall recall
the main notions of the algebraic integration theory. Then Section \ref{sec:weak-cvgce}
will be devoted to the weak convergence, divided into the tightness result
(Section \ref{sec:tightness}) and the fdd convergence (Section \ref{sec:fdd-cvgce}).
Finally, Section \ref{sec:tech} contains the technical lemmas of the paper.


\section{Background on algebraic integration and fractional SDEs}\label{sec:basic-alg-intg}

This section contains a summary of the algebraic integration introduced in
\cite{Gu}, which was also used in \cite{NNRT,NNT} in order to solve and analyze
fractional SDEs. We recall its main features here, since our approximation result
will also be obtained in this setting.

\vspace{0.2cm}

Let $x$ be a H\"older continuous $\R^d$-valued function of order  $\gamma$,
with
$1/3<\gamma\le 1/2$,
and $\sigma:\R^n\rightarrow \R^{n\times d}$, $b:\R^n\rightarrow\R^n$ be two
bounded  and smooth functions. We shall consider in the sequel the $n$-dimensional equation
\begin{equation}\label{eq1}
dy_t=\si\lp y_{t} \rp dx_t
+ b\lp   y_{t} \rp dt,\quad y_0=a\in\R^n,\quad t\in[0,T].
\end{equation}
In order to define rigorously and solve this equation,
we will need some algebraic and analytic notions which are introduced in the next
subsection.

\vspace{0.2cm}

\subsection{Increments}\label{incr}

We first present the basic  algebraic structures which
will allow us to define a pathwise integral with respect to
irregular functions. For an arbitrary real number
$T>0$, a vector space $V$ and an integer $k\ge 1$ we denote by
$\cac_k(V)$ the set of functions $g : [0,T]^{k} \to V$ such
that $g_{t_1 \cdots t_{k}} = 0$
whenever $t_i = t_{i+1}$ for some $i\le k-1$.
Such a function will be called a
\emph{$(k-1)$-increment}, and we will
set $\cac_*(V)=\cup_{k\ge 1}\cac_k(V)$. An important elementary operator
is defined by
\begin{equation}
  \label{eq:coboundary}
\delta : \cac_k(V) \to \cac_{k+1}(V), \qquad
(\delta g)_{t_1 \cdots t_{k+1}} = \sum_{i=1}^{k+1} (-1)^{k-i}
g_{t_1  \cdots \hat t_i \cdots t_{k+1}} ,
\end{equation}
where $\hat t_i$ means that this particular argument is omitted.
A fundamental property of $\der$, which is easily verified,
is that
$\delta \delta = 0$, where $\delta \delta$ is considered as an operator
from $\cac_k(V)$ to $\cac_{k+2}(V)$.
 We will denote $\cz\cac_k(V) = \cac_k(V) \cap \text{Ker}\delta$
and $\cb \cac_k(V) =
\cac_k(V) \cap \text{Im}\delta$.

\vspace{0.2cm}

Some simple examples of actions of $\der$
 are obtained  for
$g\in\cac_1(V)$ and $h\in\cac_2(V)$. Then, for any $s,u,t\in\ott$, we have
\begin{equation}
\label{eq:simple_application}
  (\der g)_{st} = g_t - g_s,
\quad\mbox{ and }\quad
(\der h)_{sut} = h_{st}-h_{su}-h_{ut}.
\end{equation}
Furthermore, it is easily checked that
$\cz \cac_{k+1}(V) = \cb \cac_{k}(V)$ for any $k\ge 1$.
In particular, the following basic property holds:
\begin{lemma}\label{exd}
Let $k\ge 1$ and $h\in \cz\cac_{k+1}(V)$. Then there exists a (non unique)
$f\in\cac_{k}(V)$ such that $h=\der f$.
\end{lemma}

Observe that Lemma \ref{exd} implies that all  elements
$h \in\cac_2(V)$  with $\der h= 0$ can be written as $h = \der f$
for some (non unique) $f \in \cac_1(V)$. Thus we get a heuristic
interpretation of $\der |_{\cac_2(V)}$:  it measures how much a
given 1-increment  is far from being an  exact increment of a
function, i.e., a finite difference.

\vspace{0.2cm}

 Note that our further discussion will mainly rely on
$k$-increments with $k \leq 2$.
For the simplicity of the exposition, we will assume from now that $V=\R^d$.
We measure the size of these increments by H\"older norms,  which are
defined in the following way: for $f \in \cac_2(V)$ let
$$
\norm{f}_{\mu} =
\sup_{s,t\in\ott}\frac{|f_{st}|}{|t-s|^\mu},
\quad\mbox{and}\quad
\cac_2^\mu(V)=\lcl f \in \cac_2(V);\, \norm{f}_{\mu}<\infty  \rcl.
$$
Obviously, the usual H\"older spaces $\cac_1^\mu(V)$  are determined
        in the following way: for a continuous function $g\in\cac_1(V)$ simply set
\begin{equation}\label{def:hnorm-c1}
\|g\|_{\mu}=\|\der g\|_{\mu},
\end{equation}
and we will say that $g\in\cac_1^\mu(V)$ iff $\|g\|_{\mu}$ is finite.
Note that $\|\cdot\|_{\mu}$ is only a semi-norm on $\cac_1(V)$,
but we will  work  in general on spaces of the type
\begin{equation}\label{def:hold-init}
\cac_{1,a}^\mu(V)=
\lcl g:\ott\to V;\, g_0=a,\, \|g\|_{\mu}<\infty \rcl,
\end{equation}
for a given $a\in V$, on which $\|g\|_{\mu}$  is a norm.
 For $h \in \cac_3(V)$ set in the same way
\begin{eqnarray}
  \label{eq:normOCC2}
  \norm{h}_{\gamma,\rho} &=& \sup_{s,u,t\in\ott}
\frac{|h_{sut}|}{|u-s|^\gamma |t-u|^\rho}\\
\|h\|_\mu &= &
\inf\left \{\sum_i \|h_i\|_{\rho_i,\mu-\rho_i} ;\, h =
 \sum_i h_i,\, 0 < \rho_i < \mu \right\} ,\nonumber
\end{eqnarray}
where the  infimum is taken over all sequences $\{h_i \in \cac_3(V) \}$
such that $h
= \sum_i h_i$ and for all choices of the numbers $\rho_i \in
(0,\mu)$.
Then  $\|\cdot\|_\mu$ is easily seen to be a norm on $\cac_3(V)$, and we set
$$
\cac_3^\mu(V):=\lcl h\in\cac_3(V);\, \|h\|_\mu<\infty \rcl.
$$
Eventually,
let $\cac_3^{1+}(V) = \cup_{\mu > 1} \cac_3^\mu(V)$,
and  note that the same kind of norms can be considered on the
spaces $\cz \cac_3(V)$, leading to the definition of  the  spaces
$\cz \cac_3^\mu(V)$ and $\cz \cac_3^{1+}(V)$.

\vspace{0.2cm}

With these notations in mind,
the crucial point in the current approach to pathwise integration of irregular
paths is that the operator
$\delta$ can be inverted under mild  smoothness assumptions. This
inverse is called $\laa$. The proof of the following proposition  may be found in
\cite{Gu}, and in a more elementary form in \cite{GT}:
\begin{proposition}
\label{prop:Lambda}
There exists a unique linear map $\Lambda: \cz \cac^{1+}_3(V)
\to \cac_2^{1+}(V)$ such that
$$
\delta \Lambda  = \id_{\cz \cac_3^{1+}(V)}
\quad \mbox{ and } \quad \quad
\Lambda  \delta= \id_{\cac_2^{1+}(V)}.
$$
In other words, for any $h\in\cac^{1+}_3(V)$ such that $\der h=0$
there exists a unique $g=\laa(h)\in\cac_2^{1+}(V)$ such that $\der g=h$.
Furthermore, for any $\mu > 1$,
the map $\laa$ is continuous from $\cz \cac^{\mu}_3(V)$
to $\cac_2^{\mu}(V)$ and we have
\begin{equation}\label{ineqla}
\|\Lambda h\|_{\mu} \leq \frac{1}{2^\mu-2} \|h\|_{\mu} ,\qquad h \in
\cz \cac^{\mu}_3(V).
\end{equation}
\end{proposition}
Moreover, $\laa$  has a nice interpretation in terms of generalized Young integrals:
\begin{corollary}
\label{cor:integration}
For any 1-increment $g\in\cac_2 (V)$ such that $\der g\in\cac_3^{1+}(V)$
set
$
\delta f = (\id-\Lambda \delta) g
$.
Then
$$
(\delta f)_{st} = \lim_{|\Pi_{ts}| \to 0} \sum_{i=0}^n g_{t_i\,t_{i+1}},
$$
where the limit is over any partition $\Pi_{st} = \{t_0=s,\dots,
t_n=t\}$ of $[s,t]$, whose mesh tends to zero. Thus, the
1-increment $\delta f$ is the indefinite integral of the 1-increment $g$.
\end{corollary}


\vspace{0.2cm}

\subsection{Weakly controlled paths}\label{sec:cl4}

This subsection is devoted to the definition of generalized integrals
with respect to a rough path of order 2, and to the resolution of
equation (\ref{eq1}). Notice that,
in the sequel of our paper, we will use  both the notations $\ist f dg$
or $\cj_{st}(f\, dg)$ for the integral of a function $f$ with
respect to a given increment $dg$ on the interval $[s,t]$.
 The second notation $\cj_{st}(f\, dg)$ will be used to avoid
some cumbersome notations in  our computations.
Observe also that the drift
term $b$ is generally harmless if one  wants to
solve the equation (\ref{eq1}). See e.g. Remark 3.14 in \cite{NNRT}.
 Hence, we will simply deal with
an equation of the form
\begin{equation}\label{red:eq}
dy_t=\si\lp y_{t} \rp dx_t,
\quad t\in\ott,
\quad\mbox{ with }\quad
y_0=a
\end{equation}
in the remainder of this section.
\vspace{0.2cm}

Before going into the technical  details,
let us make some heuristic considerations about the properties that
a  solution  of equation   (\ref{eq1}) should have. Set
$\hsi_t=\si\lp y_{t} \rp$, and suppose that $y$ is a solution of  (\ref{red:eq}),
with $y\in\cac_1^\ka$ for a given $1/3<\ka<\ga$.
Then the integral form of our equation can be
written as
\begin{equation}\label{intg}
y_t=a+\iot \hsi_u dx_u,  \qquad t\in\ott.
\end{equation}
Our approach to generalized integrals induces us to work with increments
of the form $(\delta y)_{st}=y_t-y_s$  instead  of (\ref{intg}). However, it is easily checked that one can decompose
(\ref{intg}) into
$$
(\delta y)_{st}=\ist \hsi_u dx_u=\hsi_s (\delta x)_{st}+\rho_{st},
\quad\mbox{ with }\quad
\rho_{st}=\ist (\hsi_u-\hsi_s) dx_u,
$$
 if our  integral is linear.
We thus have  obtained a decomposition of $y$ of the form
$\delta y=\hsi\delta x+\rho$.  Let us see, still at a heuristic level,
 which regularity we can expect for $\hsi$ and $r$. If
$\si$ is a $C_b^1$-function, we have that $\hsi$ is bounded and
$$
|\hsi_t-\hsi_s|
\leq
\|\nabla\si  \|_{\infty}  \|y\|_{\ka} |t-s|^{\ka},
$$
where $\|y\|_{\ka}$ denotes the H\"older norm of $y$  defined by
(\ref{def:hnorm-c1}).  Hence
we have that $\hsi$ belongs to $\cac_1^\ka$ and is bounded.
As far as $\rho$ is
concerned, it should
inherit both the regularities of $\delta\hsi$ and $x$,  provided that  the integral
$\ist (\hsi_u-\hsi_s) dx_u=\ist(\delta\hsi)_{su}dx_u$ is well defined.  Thus, one should
expect that $\rho\in\cac_2^{2\ka}$, and even $\rho\in\cac_2^{\ka+\ga}$.
To summarize, we have found that
a solution $\delta y$ of the equation should be decomposable into
\begin{equation}\label{first:structure1}
\delta y =\hsi \delta x + \rho,
\quad\mbox{ with }\quad
\hsi\in\cac_1^\ga
\mbox{ bounded and }
\rho\in\cac_2^{2\ka}.
\end{equation}
This is precisely the structure we will  demand for a possible
solution  of
(\ref{red:eq}):
\begin{definition}\label{def:ccp}
Let $z$ be a path in $\cac_1^\ka(\R^k)$ with $\ka\le\ga$
and $2\ka+\ga>1$.
We say that $z$ is a controlled path based on $x$, if
$z_0=a$, which is a given initial condition in $\R^k$,
and $\der z\in\cac_2^\ka(\R^k)$ can be decomposed into
\begin{equation}\label{weak:dcp}
\der z=\zeta \der x+ r,
\quad\mbox{i.\!\! e.}\quad
(\der z)_{st}=\zeta_s (\der x)_{st} + \rho_{st},
\quad s,t\in\ott,
\end{equation}
with $\zeta\in\cac_1^\ka(\R^{k\times d})$ and $\rho$ is a regular part
belonging to $\cac_2^{2\ka}(\R^k)$.
The space of controlled
paths will be denoted by $\cq_{\ka,a}(\R^k)$, and a path
$z\in\cq_{\ka,a}(\R^k)$ should be considered in fact as a couple
$(z,\zeta)$. The natural semi-norm on $\cq_{\ka,a}(\R^k)$ is given
by
$$
\cn[z;\cq_{\ka,a}(\R^k)]=
\cn[z;\cac_1^{\ka}(\R^k)]
+ \cn[\zeta;\cac_1^{b}(\R^{k,d})]
+ \cn[\zeta;\cac_1^{\ka}(\R^{k,d})]
+\cn[\rho;\cac_2^{2\ka}(\R^k)]
$$
with $\cn[g;\cac_1^{\ka}(V)]$ defined by
(\ref{def:hnorm-c1}) and
$\cn[\zeta;\cac_1^{b}(V)]=\sup_{0\le s\le T}|\zeta_s|_V$.
\end{definition}

\vspace{0.2cm}

Having defined our algebraic and analytic framework, we now can give a
sketch of the strategy used in \cite{Gu} in order to solve
equation (\ref{red:eq}):
\begin{enumerate}
\item
Verify the stability of $\cq_{\ka,a}(\R^k)$ under a smooth map
$\vp:\R^k\to\R^n$.
\item
Define rigorously the integral $\int z_u dx_u=\cj(z dx)$
for a controlled path $z$ and computed its decomposition
(\ref{weak:dcp}).
\item
Solve equation (\ref{red:eq}) in the space $\cq_{\ka,a}(\R^k)$
by a fixed point argument.
\end{enumerate}
Actually, for the second point  one has
to assume a  priori the following hypothesis  on the driving rough
  path, which is standard  in rough path type considerations:
\begin{hypothesis}\label{hyp:x-cl}
The $\R^d$-valued
$\ga$-H\"older path $x$ admits a L\'evy area,
 that is a process $\xd=\cj(dx dx)\in\cac_2^{2\ga}(\R^{d\times d})$
satisfying
$$
\der\xd=\der x\otimes \der x,
\quad\mbox{i.\!\! e.}\quad
\lc (\der\xd)_{sut} \rc(i,j)
=
[\der x^{i}]_{su} [\der x^{j}]_{ut},
\quad s,u,t\in\ott, \, i,j\in\{1,\ldots,d  \}.
$$
\end{hypothesis}
Then the following result is proved in \cite{Gu}, using the strategy sketched above:
\begin{theorem}\label{thm:ex-uniq1}
Let $x$ be a process satisfying Hypothesis \ref{hyp:x-cl} and
$\si:\R^n\to\R^{n\times d}$ be a $C^2$ function, which is bounded
together with its derivatives. Then
\begin{enumerate}
\item
Equation (\ref{red:eq}) admits a unique solution $y$ in
$\cq_{\ka,a}(\R^n)$ for any $\ka<\ga$ such that $2\ka+\ga>1$.
\item
The mapping $(a,x,\xd)\mapsto y$ is continuous from
$\R^n\times\cac_1^{\ga}(\R^d)\times\cac_2^{2\ga}(\R^{d\times d})$
to $\cq_{\ka,a}(\R^n)$.
\end{enumerate}
\end{theorem}
We shall see in the next subsection that this general theorem can be applied in
the fBm context.

\vspace{0.2cm}

\subsection{Application to the fBm}\label{sec:mall-calculus}

Let $B=(B^1,\ldots,B^d)$ be a $d$-dimensional Liouville fBm of Hurst
index $H\in(\frac13,\frac12)$, that is $B^1,\ldots,B^d$ are
$d$ independent centered Gaussian processes of the form
$$
B_t^i=\iot (t-r)^{H-\frac12} dW_r^i,
$$
where $W=(W^1,\ldots,W^d)$ is a $d$-dimensional Wiener process.
The next lemma will be useful all along the paper.
\begin{lemma}\label{onenabesoin}
There exists a positive constant $c$, depending only on $H$, such that
\begin{equation}
E|B^i_t-B^i_s|^2=\int_0^s \big[(t-r)^{H-\frac12}-(s-r)^{H-\frac12}\big]^2 dr
+ \int_s^t (t-r)^{2H-1} dr
 \leq c|t-s|^{2H}\label{esti1}
\end{equation}
for all $t>s\geq 0$.
\end{lemma}
\begin{proof}
Indeed, it suffices to observe that 
\begin{eqnarray*}
\int_0^s \big[(t-r)^{H-\frac12}-(s-r)^{H-\frac12}\big]^2 dr
&=&
\int_0^s \big[(t-s+r)^{H-\frac12}-r^{H-\frac12}\big]^2 dr \\
&=& (t-s)^{2H}
\int_0^{\frac{s}{t-s}} \big[(1+r)^{H-\frac12}-r^{H-\frac12}\big]^2 dr\\
&\leq& (t-s)^{2H}
\int_0^{\infty} \big[(1+r)^{H-\frac12}-r^{H-\frac12}\big]^2 dr
\end{eqnarray*}
and that $\int_s^t (t-r)^{2H-1}dr=\frac{(t-s)^{2H}}{2H}$.

\end{proof}

Let $\ce$ be the set of step-functions on $[0,T]$ with values in $\R^{d}$.  Consider the Hilbert space $\ch$ defined as the closure of $\ce$ with respect to the scalar product induced by
$$
{ \left \lpa (\1_{[0,t_{1}]}, \ldots, \1_{[0,t_{d}]}), (\1_{[0,s_{1}]}, \ldots, \1_{[0,s_{d}]}) \right\rpa}_{\ch}\; =\; \sum_{i=1}^{d}R(t_{i},s_{i}), \quad s_{i},t_{i} \in [0,T], \,\, i=1, \ldots, d,
$$
where $R(t,s):=E[B^i_tB^i_s]$.
Then a natural representation of the inner product in $\ch$ is given via the operator
$\mathscr{K}$, defined from $\ce$ to $L^2(\ott)$, by:
$$
\mathscr{K}\vp(t)=(T-t)^{H-\frac12}\vp(t) - \lp \frac12-H\rp\int_t^T [\vp(r)-\vp(t)] (r-t)^{H-\frac32} \, dr,
$$
and it can be checked that $\mathscr{K}$ can be extended as an isometry between $\ch$
and the Hilbert space $L^2(\ott;\R^d)$. Thus the inner product in $\ch$ can be defined as:
$$
\lla\vp,\psi\rra_{\ch} \triangleq \lla \mathscr{K}\vp,\mathscr{K}\psi\rra_{L^2(\ott;\R^d)}.
$$
The mapping
$
(\1_{[0,t_{1}]}, \ldots, \1_{[0,t_{d}]}) \mapsto \sum_{i=1}^{d}B_{t_{i}}^{i}
$
can also be extended into an isometry between $\ch$ and the first
Gaussian chaos $H_{1}(B)$ associated with $B=(B^{1}, \ldots , B^{d})$.
We denote this isometry by $\varphi \mapsto B(\varphi)$, and $B(\varphi)$ is called
the Wiener-It\^o integral of $\varphi$. It is shown in \cite[page 284]{fdp} that
$\cac_1^\ga(\R^d)\subset\ch$ whenever
$\ga>1/2-H$, which allows to define $B(\vp)$ for such kind of functions.

\vspace{0.2cm}

We are now ready to prove that Theorem \ref{thm:ex-uniq1} can be applied to the
Liouville fBm, which amounts to check Hypothesis \ref{hyp:x-cl}.
\begin{proposition}\label{prop:hyp-fbm}
Let $B$ be a $d$-dimensional Liouville fBm,
and suppose that its Hurst parameter satisfies $H\in(1/3,1/2)$. Then almost all sample paths of  $B$ satisfy
Hypothesis \ref{hyp:x-cl}, with any Hölder exponent $1/3<\ga<H$, and
a Lévy area given by
$$
\bdst=\ist dB_u \otimes \int_s^u dB_v,
\quad\mbox{i. e.}\quad
\bdst(i,j)=\ist dB_u^i  \int_s^u dB_v^j,
\quad i,j\in\{1,\ldots,d  \},
$$ for $0\le s < t \leq T$. Here,
the stochastic integrals are defined as Wiener-It\^o integrals when $i\neq j$, while, when $i=j$, they are
simply given by
$$
\ist dB_u^i  \int_s^u dB_v^i = \frac12\left(B_t^i - B_s^i\right)^2.
$$
\end{proposition}

\begin{proof}
First of all, it is a classical fact
that $B\in\cac_1^\ga(\R^d)$ for any $1/3<\ga<H$, when $B$ is a Liouville fBm with $H>1/3$
(indeed, combine the Kolmogorov-$\check{\rm C}$entsov theorem with Lemma \ref{onenabesoin}). 
Furthermore, we have already mentioned that $\cac_1^\ga(\R^d)\subset\ch$ for any $\ga>1/2-H$.
In particular, if $H>\ga>1/3$, the condition  $\ga>1/2-H$ is satisfied and, conditionaly to $B^j$, 
$\ist dB_u^i  \int_s^u dB_v^j$ is well-defined for $i\neq j$, as a Wiener-It\^o integral with
respect to $B^i$, of the form $B^i(\vp)$ for a well-chosen $\vp$. Hence, $\bd$
is almost surely a well-defined element of $\cac_2(\R^{d\times d})$.

\vspace{0.2cm}

Now, simple algebraic computations immediately
yield  that $\der\bd=\der B\otimes\der B$.
Furthermore, Lemma \ref{lm-technical} yields
$$
E\lc |\bdst(i,j)|^2  \rc \leq c |t-s|^{4H}.
$$
Invoking this inequality and thanks to the fact that $\bd$ is a process in the
second chaos of $B$, on which all
$L^p$ norms ($p>1$) are equivalent, we get that
\begin{equation*}
E\lc |\bdst(i,j)|^p  \rc
\leq c_p |t-s|^{2pH}.
\end{equation*}
This allows to conclude, thanks to an elaboration of Garsia's lemma which can be found in \cite[Lemma 4]{Gu} (and will be recalled at (\ref{lem:garsia})),
that $\bd\in\cac_2^{2\ga}(\R^{d\times d})$ for any
$\ga<1/3$. This ends the proof.

\end{proof}

With all these results in hand, we have obtained a reasonable definition of diffusion processes driven by a fBm, and we can now proceed to their approximation in law.


\section{Approximating sequence}\label{sec:weak-cvgce}

In this section, we will introduce our smooth approximation of $B$, namely $X^\e$, which
shall converge {\it in law} to $B$. This will allow to interpret equation (\ref{eq:eds-approx}) in
the usual Lebesgue-Stieltjes sense. We will then study the convergence in law of the process
$y^\e$ solution to (\ref{eq:eds-approx}) towards the solution $y$ of (\ref{eq:eds-intro}).

\vspace{0.2cm}

As mentioned in the introduction, the approximation of $B$ we shall deal with is defined 
as follows, for $i=1,\ldots,d$:
\begin{equation}\label{eq:def2-x-ep}
X^{i,\e}(t)=\int_0^t (t+\e-r)^{H-\frac12}\theta^{\e,i}(r)dr,
\end{equation}
where
$$
\theta^{\e,i}(r)=\frac1{\e}(-1)^{N^i(\frac{r}{\e})},
$$
for $N^i$, $i=1,\ldots, d$, some independent standard Poisson processes.
Furthermore, we have recalled in Theorem \ref{thm:ex-uniq1} that the solution
$y$ to (\ref{eq:eds-intro}) is a continuous function of $(a,B,\bd)$, considered
respectively as elements of $\R^d,\cac_1^\ga(\R^d)$ and $\cac_2^{2\ga}(\R^{d\times d})$ for
$1/3<\ga<H$. Thus our approximation theorem \ref{thm:approx-diffusion}
can be easily deduced from the following result:
\begin{theorem}\label{main-thm-weak}
For any $\e>0$, let ${\bf X}^{\bf 2,\e}=({\bf X}^{\bf 2,\e}_{st}(i,j))_{s,t\ge 0;\,i,j=1,\ldots,d}$
be the natural L\'evy's area associated to $X^\e$,
defined by
\begin{equation}\label{levydef}
{\bf X}^{\bf 2,\e}_{st}(i,j)=\int_s^t (X^{j,\e}_u-X^{j,\e}_s)dX^{i,\e}_u,
\end{equation}
where the integral is understood in the usual Lebesgue-Stieltjes sense.
Then, as $\e\rightarrow 0$,
\begin{equation}\label{loi}
(X^\e,{\bf X}^{\bf 2,\e})
\,\,\,{\stackrel{{\rm Law}}{\longrightarrow}}\,\,\,
(B,{\bf B}^{\bf 2}),
\end{equation}
where ${\bf B}^{\bf 2}$ denotes the L\'evy area defined in Proposition \ref{prop:hyp-fbm},
and where the convergence in law holds in spaces
$\mathcal{C}_1^\mu(\R^d)\times\mathcal{C}_2^{2\mu}(\R^{d\times d})$,
for any $\mu<H$.
\end{theorem}

The remainder of our work is devoted to the proof of Theorem \ref{main-thm-weak}.
As usual in the context of weak convergence of stochastic processes, we divide the proof into the 
weak convergence for finite-dimensional distributions (Section \ref{sec:fdd-cvgce}) and
a tightness type result (Section \ref{sec:tightness}). 

\begin{remark}
A natural idea for the proof of Theorem \ref{main-thm-weak} could be 
to use the methodology initiated by Kurtz and Protter in \cite{KP}.
But the problem, here, is that the quantities we are dealing with are
not ``close enough" to a martingale.
\end{remark}

\section{Tightness in Theorem \ref{main-thm-weak}}\label{sec:tightness}

From now, we write $\cac_1^\mu$ (resp. $\cac_2^{2\mu}$) instead of $\cac_1^\mu(\R^d)$ (resp. $\cac_2^{2\mu}(\R^{d\times d})$). 
We first need a general tightness criterion in the Hölder spaces $\cac_1^\mu$ and $\cac_2^{2\mu}$.
\begin{lemma}\label{lem:compactness-in-c12}
Let $\mathscr{E}^\gamma$ denote the set of $(x,{\bf x}^{\bf 2})
\in
\mathcal{C}_1^\gamma\times\mathcal{C}_2^{2\gamma}$ verifying $x_0=0$ and 
\begin{equation}\label{are}
\forall s,t\ge 0,\,\forall i,j=1,\ldots,d:\quad
{\bf x}^{\bf 2}_{st}(i,j)
={\bf x}^{\bf 2}_{0t}(i,j)-{\bf x}^{\bf 2}_{0s}(i,j)-x_s^i(x^j_t-x^j_s).
\end{equation}
Let $\mu$ such that $0\le \mu <\gamma$. Then, any bounded subset $\mathscr{Q}$
of $\mathscr{E}^\gamma$ is precompact in $\mathcal{C}_1^\mu\times\mathcal{C}_2^{2\mu}$.
\end{lemma}
\begin{proof} Let $(x^n,{\bf x^{2,{\it n}}})$ be a sequence of $\mathscr{Q}$.
By assumption, $(x^n,{\bf x}^{\bf 2,{\it n}}_{0\cdot})$ is bounded and equicontinuous.
Then, Ascoli's theorem applies and, at least along a subsequence, which may also
be called
$(x^n,{\bf x}^{\bf 2,{\it n}}_{0\cdot})$, it converges uniformly to $(x,{\bf x}^{\bf 2}_{0\cdot})$.
Using (\ref{are}), we obtain in fact that $(x^n,{\bf x}^{\bf 2,{\it n}})$ converges
uniformly to $(x,{\bf x}^{\bf 2})$. Moreover, since we obviously have
$$
\|x\|_\mu\le\liminf_{n\rightarrow\infty} \|x^n\|_\mu\quad\mbox{and}\quad
\|{\bf x^2}\|_{2\mu}\leq\liminf_{n\rightarrow\infty} \|{\bf x}^{\bf 2,{\it n}}\|_{2\mu},
$$
we deduce that $(x,{\bf x}^{\bf 2})\in\mathcal{C}_1^\mu\times\mathcal{C}_2^{2\mu}$.
Finally, we have
$$
\|x-x^n\|_\mu\longrightarrow 0\quad\mbox{and}\quad\|{\bf x^2}-{\bf x^{2,{\it n}}}\|_{2\mu}
\longrightarrow 0,
$$
owing to the fact that
$$
\|x-x^n\|_\mu\le\|x-x^n\|_\gamma \,{\|x-x^n\|_\infty}^{1-\frac\mu\gamma}
\leq \big(\|x\|_\gamma + \|x^n\|_\gamma\big){\|x-x^n\|_\infty}^{1-\frac\mu\gamma}
$$
and similarly:
$$
\|{\bf x^2}-{\bf x^{2,{\it n}}}\|_{2\mu}
\leq \big(\|{\bf x^2}\|_{2\gamma} + \|{\bf x^{2,{\it n}}}\|_{2\gamma}\big)
{\|{\bf x^2}-{\bf x^{2,{\it n}}}\|_\infty}^{1-\frac\mu\gamma}.
$$
\end{proof}

We will use the last result in order to get a reasonable tightness criterion
for our approximation processes $X^\e$ and ${\bf X^{2,\e}}$, by means
of a slight elaboration of \cite[Corollary 6.1]{lejay}:
\begin{proposition}\label{tight}
Let $X^\e$ and ${\bf X^{2,\e}}$ be defined respectively by (\ref{eq:def2-x-ep}) and (\ref{levydef}).
If, for every $\eta>0$, there exists $\gamma>\mu$ and $A<\infty$ such that
\begin{equation}\label{tight-cond}
\sup_{0<\e\le 1} P[\|X^\e\|_\gamma>A]\leq\eta\quad\mbox{and}\quad
\sup_{0<\e\le 1} P[\|{\bf X^{2,\e}}\|_{2\gamma}>A]\leq\eta,
\end{equation}
then $(X^\e,{\bf X^{2,\e}})$ is tight in
$\mathcal{C}_1^\mu\times\mathcal{C}_2^{2\mu}$.
\end{proposition}
\begin{proof}
Recall the Prokhorov theorem relating precompactness of measures {\it on} a space
to compactness of sets {\it in} the space. This result states that a family $M$
of probability measures on the Borel sets of a complete separable metric space $S$
is weakly precompact if and only if for every $\eta>0$ there exists a compact set
$K_\eta\subset S$ such that
$$
\sup_{\mu\in M} \mu\lp S\setminus K_\eta\rp\leq \eta.
$$
Furthermore, it is readily checked that the couple $(X^\e,{\bf X^{2,\e}})$ satisfies
the assumption (\ref{are}), which allows to apply Lemma \ref{lem:compactness-in-c12}.
Hence, combining this lemma with Prokhorov's theorem, our proposition is easily proved.

\end{proof}

\vspace{0.2cm}

Let us turn now to the main result of this subsection:
\begin{proposition}
The sequence $(X^\e,{\bf X}^{\bf 2,\e})_{\ep>0}$ defined in Theorem \ref{main-thm-weak} is tight in
$\mathcal{C}_1^\mu\times\mathcal{C}_2^{2\mu}$.
\end{proposition}

\begin{proof}
Thanks to Proposition \ref{tight}, we just have to
prove that $(X^\e,{\bf X}^{\bf 2,\e})$ verifies (\ref{tight-cond}).
For an arbitrary $\eta\in(0,1)$, we will first deal with the relation
\begin{equation}\label{eq:bnd-hold-norm-xe}
\sup_{0<\e\le 1} P\big[\|X^\e\|_\gamma>A\big]\leq\eta,
\end{equation}
for $A=A_\eta$ large enough, and $1/3<\ga<H$. To this purpose, let us recall some
basic facts about Sobolev spaces, for which we refer to \cite{Ad}: for
$\al\in(0,1)$ and $p\ge 1$, the Sobolev space $\cw^{\al,p}(\ott;\,\R^n)$ is induced by the
semi-norm
\begin{equation}\label{eq:def-sobolev-norm}
\|f\|_{\al,p}^p=\iott\iott \frac{|f(t)-f(s)|^p}{|t-s|^{1+\al p}} \, ds dt.
\end{equation}
Then the Sobolev imbedding theorem states that, if $\al p>1$, then $\cw^{\al,p}(\ott;\,\R^d)$
is continuously imbedded in $\cac_1^{\ga}(\R^d)$ for any $\ga<\al-1/p$, where
the spaces $\cac_1^{\ga}$ have been defined by relation (\ref{def:hnorm-c1}), and in this case, we furthermore have that
\begin{equation}\label{eq:24}
\|f\|_{\ga} \le c \|f\|_{\al,p},
\end{equation}
for a positive constant $c=c_{\al,p}$. Notice that, in both (\ref{def:hnorm-c1}) and (\ref{eq:def-sobolev-norm}), the $\sup$ part
of the usual Hölder or Sobolev norm has been omitted, but can be recovered since
we are dealing with fixed initial conditions.
In order to prove (\ref{eq:bnd-hold-norm-xe}), it is thus sufficient to check that,
for any $p \geq 1$ sufficiently large and $\al<H$, the following bound holds true:
\begin{equation}\label{eq:bnd-sob-norm-xe}
\sup_{0<\e\le 1} E \lc \iott\iott \frac{|X^\e(t)-X^\e(s)|^p}{|t-s|^{1+\al p}} \, ds dt \rc
\leq M_{\al,p} <\infty.
\end{equation}
Invoking Lemma \ref{lm-control}, we get, for any $\e>0$, any $t>s\ge 0$ and any integer $m\geq 1$:
\begin{eqnarray}
&&E\left[|X^{\e,i}(t)-X^{\e,i}(s)|^{2m}\right]\\
&\leq& 2^{2m-1}
E\left[\left|\int_0^s \big( (t+\e-r)^{H-\frac12}-(s+\e-r)^{H-\frac12}\big)\theta^{\e,i}(r)dr\right|^{2m}\right]\notag\\
&&+ 2^{2m-1}
E\left[\left|\int_s^t (t+\e-r)^{H-\frac12}\,\theta^{\e,i}(r)dr\right|^{2m}\right]\notag\\
&\leq& \frac{2^{m-1}(2m)!}{m!}\left( \int_0^s \big( (t+\e-r)^{H-\frac12}-(s+\e-r)^{H-\frac12}\big)^2dr
\right)^{m}\notag\\
&&+ \frac{2^{m-1}(2m)!}{m!}\left( \int_s^t (t+\e-r)^{2H-1}dr\right)^{m}\notag\\
&\leq& \frac{2^{m-1}(2m)!}{m!}\left( \int_0^s \big( (t-r)^{H-\frac12}-(s-r)^{H-\frac12}\big)^2dr
\right)^{m}\\
&&+ \frac{2^{m-1}(2m)!}{m!}\left( \int_s^t (t-r)^{2H-1}dr\right)^{m}\notag\\
&\leq & c_{2m,H}|t-s|^{2mH} \label{cont}\quad\mbox{by Lemma \ref{onenabesoin}}.
\end{eqnarray}
Note that here, and in the remainder of the proof, 
$c_{\{\cdot\}}$
denotes a generic constant depending only on the object(s) inside its argument, and which may take different values one
formula to another one.
From (\ref{cont}), we deduce that
(\ref{eq:bnd-sob-norm-xe}) holds for any $\al<H$ and $p$ large enough, from which
(\ref{eq:bnd-hold-norm-xe}) is easily seen.
Moreover, thanks to the classical Garsia-Rodemich-Rumsey lemma, see \cite{GRR}, for any $\e,\delta,T>0$ 
and $i\in\{1,\ldots,d\}$,
there exists a random variable $G^{T,\delta,\e,i}$
such that, for any $s,t\in [0,T]$:
\begin{equation}\label{GRR}
|X^{\e,i}(t)-X^{\e,i}(s)|\leq G^{T,\delta,\e,i} |t-s|^{H-\delta}.
\end{equation}
Since the bound in (\ref{cont}) is {\sl independent} of $\e$, it is easily checked
that, for any integer $m\geq 1$, any $i\in\{1,\ldots,d\}$ and any $\delta,T>0$
($\delta$ small enough), we have
$$c_{2m,\delta}:=\sup_{0<\e\leq 1}E\left(\vert G^{T,\delta,\e,i}\vert^{2m}\right)< +\infty.$$

\smallskip

Let us turn now to the tightness of $({\bf X}^{\bf 2,\e})_{\ep>0}$. Recall first that
$
{\bf X}^{\bf 2,\e}_{st}(i,i) = \frac12( X^{\e,i}_t-X^{\e,i}_s)^2.
$
Therefore, we deduce from (\ref{cont}) that
\begin{equation}
E\big[|{\bf X}^{\bf 2,\e}_{st}(i,i)|^{2m}\big]\leq \frac{c_{4m,H}}{2^{2m}}|t-s|^{4mH}.\label{ineq:norm-xd1}
\end{equation}
Assume now that $i\neq j$. We have, by applying successively (\ref{lm2}), 
Lemma \ref{lm-control} and (\ref{GRR}):
\begin{eqnarray}
E[|{\bf X}^{\bf 2,\e}_{st}(i,j)|^{2m}]
&\leq &c_m\, E\left|\int_s^t  \big(X_u^{j,\e}-X_s^{j,\e}\big)^2 (t+\e-u)^{2H-1} du
\right|^m \notag\\
&&+c_m\,E\left|\int_0^s  \big(X_u^{j,\e}-X_s^{j,\e}\big)^2 \big((t+\e-u)^{H-\frac12}-(s+\e-u)^{H-\frac12}\big)^2du\right|^m \notag\\
&&+c_{m,H}\,E\left|\int_0^t \left( \int_{s\vee v}^t 
\vert X_u^{j,\e}-X_v^{j,\e}\vert\,\, (u+\e-v)^{H-\frac32} du\right)^2 dv\right|^m\notag
\end{eqnarray}
This last expression can be trivially bounded by considering the case $\e=0$, and some elementary calculations then lead to the relation
\begin{equation}\label{ineq:norm-xd2}
E[|{\bf X}^{\bf 2,\e}_{st}(i,j)|^{2m}] \le
c_{m,H}\,|t-s|^{4mH-2m\delta}.
\end{equation}
In order to conclude that ${\bf X}^2$ verifies the second inequality in (\ref{tight-cond}),
let us recall the following inequality from \cite{Gu}:
let $g\in \cac_2(V)$ for a given Banach space $V$; then,
for any $\ka>0$ and $p\ge 1$ we have
\begin{equation}\label{lem:garsia}
\| g\|_{\ka}\leq c \lp U_{\ka+2/p;p}(g) + \| \der g\|_{\ga}\rp
\quad\mbox{ with }\quad
U_{\ga;p}(g)=
\lp \int_0^T\int_0^T \frac{|g_{st}|^p}{|t-s|^{\ga p}}dsdt \rp^{1/p}.
\end{equation}
By plugging inequality (\ref{ineq:norm-xd1})-(\ref{ineq:norm-xd2}), for $\delta>0$ small enough, into (\ref{lem:garsia}) and
by recalling that $\der{\bf X}^{2,\e}=\der X^\e\otimes\der X^\e$ and inequality
(\ref{GRR}), we obtain easily the second part of (\ref{tight-cond}).

\end{proof}

\section{Fdd convergence in Theorem \ref{main-thm-weak}}\label{sec:fdd-cvgce}

This section is devoted to the second part of the proof of
Theorem \ref{main-thm-weak}, namely the convergence of finite
dimensional distributions. Precisely, we shall prove the following:
\begin{proposition}\label{prop-weak}
Let $(X^\e,{\bf X}^{\bf 2,\e})$ be the approximation process defined by
(\ref{eq:def2-x-ep}) and (\ref{levydef}). Then
\begin{equation}\label{fdd}
{\rm f.d.d.}-\lim_{\e\to 0} (X^\e,{\bf X}^{\bf 2,\e})
= (B,{\bf B}^{\bf 2}),
\end{equation}
where ${\rm f.d.d.}-\lim$ stands for the convergence in law of the
finite dimensional distributions. Otherwise stated, for any $k\ge 1$ and any family
$\{s_i,t_i;\, i\le k, 0\le s_i<t_i\le T\}$, we have
\begin{equation}\label{fddbis}
\cl-\lim_{\e\to 0}
(X_{t_1}^\e,{\bf X}_{s_1 t_1}^{\bf 2,\e},\ldots, X_{t_k}^\e,{\bf X}_{s_k t_k}^{\bf 2,\e})
= (B_{t_1},{\bf B}_{s_1 t_1}^{\bf 2},\ldots, B_{t_k},{\bf B}_{s_k t_k}^{\bf 2}).
\end{equation}
\end{proposition}
\begin{proof}
The proof is divided into several steps.

\vspace{0.2cm}

\noindent
{\it (i) Reduction of the problem}.
For simplicity, we assume that the dimension $d$ of $B$ is $2$ (the general case can be treated along the same lines,
up to some cumbersome notations). For $i=1,2$, $\e>0$ and $0\leq u\leq t\leq T$, let us consider
$$
Y^{i,\e}(u,t)=\int_u^t (X_v^{i,\e}-X_u^{i,\e})(v-u)^{H-\frac32}dv
$$
and
$$
Y^i(u,t)=\int_u^t (B^i_v - B^i_u)(v-u)^{H-\frac32}dv.
$$

In this step, we shall prove that the fdd convergence (\ref{fdd}) is a consequence of the following one:
\begin{eqnarray}
&&\left(
\int_0^\cdot \theta^{\e,1}(u)du,\int_0^\cdot \theta^{\e,2}(u)du,
\int_0^\cdot X_u^{2,\e} \theta^{\e,1}(u)du,\right.\notag\\
&&\hskip1cm\left.\int_0^\cdot Y^{2,\e}(u,\cdot)\theta^{\e,1}(u)du,
\int_0^\cdot X_u^{1,\e}\theta^{\e,2}(u)du,\int_0^\cdot Y^{1,\e}(u,\cdot)\theta^{\e,2}(u)du\right)\notag\\
&&\,\,\,{\stackrel{{\rm f.d.d.}}{\longrightarrow}}\,\,\,
\left(
W^1,W^2,\int_0^\cdot B^2_u dW_u^1,\int_0^\cdot Y^2(u,\cdot)dW^1_u,\int_0^\cdot B^1_u dW^2_u,
\int_0^\cdot Y^1(u,\cdot)dW^2_u
\right).
\label{fdd2}
\end{eqnarray}
Indeed, assume for an instant that (\ref{fdd2}) takes place. Then, approximating the kernel $(t-\cdot)^{H-1/2}$ in $L^2$ by a sequence of step functions (along the same lines as in \cite[Proof of Theorem 1, p. 404]{DJ}), it is easily checked that we also have:
\begin{eqnarray}
&&\left(
X^{1,\e},X^{2,\e},
\int_0^\cdot  (\cdot+\e-u)^{H-\frac12}\,X_u^{2,\e}\, \theta^{\e,1}(u)du,\right.\notag\\
&&\hskip1cm\left.\int_0^\cdot Y^{2,\e}(u,\cdot)\theta^{\e,1}(u)du,
\int_0^\cdot (\cdot+\e-u)^{H-\frac12}\,X_u^{1,\e}\, \theta^{\e,2}(u)du,\int_0^\cdot Y^{1,\e}(u,\cdot)\theta^{\e,2}(u)du\right)\notag\\
&&\,\,\,{\stackrel{{\rm f.d.d.}}{\longrightarrow}}\,\,\,
\left(
B^1,B^2,\int_0^\cdot (\cdot-u)^{H-\frac12}\,B^2_u\, dW_u^1,\right.\notag\\
&&\hskip2cm\left.
\int_0^\cdot Y^2(u,\cdot)dW^1_u,
\int_0^\cdot (\cdot-u)^{H-\frac12}\,B^1_u\, 
dW^2_u,\int_0^\cdot Y^1(u,\cdot)dW^2_u
\right).\notag\\
\label{fdd3}
\end{eqnarray}
In other words, we can add the deterministic kernel $(\cdot+\e-u)^{H-\frac12}$ in the first, second, 
third and fifth components of (\ref{fdd2})
without difficulty. Let us invoke now the forthcoming identity (\ref{lm2}) in Lemma \ref{otherform} for $s=0$, which allows easily to go from (\ref{fdd3}) to:
\begin{equation}\label{fdd4}
\big(X^{1,\e},X^{2,\e},{\bf X}^{\bf 2,\e}_{0\cdot}(1,2),{\bf X}^{\bf 2,\e}_{0\cdot}(2,1)\big)
\,\,\,{\stackrel{{\rm f.d.d.}}{\longrightarrow}}\,\,\,
\left(B^1,B^2,\int_0^\cdot B^2 dB^1,\int_0^\cdot B^1 dB^2\right).
\end{equation}
Finally, in order to prove our claim (\ref{fddbis}) from (\ref{fdd4}), it is enough to observe that ${\bf X}^{\bf 2,\e}_{0t}(i,i)=(X_t^{i,\e})^2/2$ and
$$
{\bf X}^{\bf 2,\e}_{st}(i,j)
={\bf X}^{\bf 2,\e}_{0t}(i,j)-{\bf X}^{\bf 2,\e}_{0s}(i,j)-X_s^{i,\e}\big(X^{j,\e}_t-X^{j,\e}_s\big).
$$

\vspace{0.2cm}

\noindent
{\it (ii) Simplification of the statement (\ref{fdd2})}. 
For the simplicity of the exposition, we only prove (\ref{fdd2}) for a fixed $t$, instead
of a vector $(t_1,\ldots,t_m)$. It will be clear from our proof that the general case can be elaborated easily from this particular situation, up  to some additional unpleasant notations. Precisely, we shall prove that, for any $u:=(u_1,\ldots,u_6)\in\R^6$, we have $\lim_{\ep\to 0}\delta_\e=E[\exp(i\langle u,\, U\rangle)]$, where $\delta_\e:=E[\exp(i\langle u,\, U^\e\rangle)]$, $U^\e$ is defined by
\begin{multline*}
U^\e=u_1\int_0^t \theta^{\e,1}(v)dv + u_2 \int_0^t \theta^{\e,2}(v)dv 
+ u_3 \int_0^t X^{2,\e}_u \theta^{\e,1}(v)dv  \\
+u_4\int_0^t Y^{2,\e}(v,t)\theta^{\e,1}(v)dv + u_5\int_0^t X^{1,\e}_v \theta^{\e,2}(v)dv
+u_6 \int_0^t Y^{1,\e}(v,t)\theta^{\e,2}(v)dv,
\end{multline*}
and 
\begin{multline*}
U=u_1W^1_t + u_2 W^2_t + u_3 \int_0^t B^2_v dW^1_v  \\
+u_4\int_0^t Y^{2}(v,t)dW^1_v + u_5\int_0^t B^1_v dW^2_v
+u_6 \int_0^t Y^{1}(v,t)dW^2_v.
\end{multline*}

\smallskip

In order to analyze the asymptotic behavior of $\delta_\e$, let us first express $U^\e$ as an integral with respect to $\theta^{\e,1}$ only. Indeed, Fubini's theorem easily yields
$$
\int_0^t X_v^{1,\e} \theta^{\e,2}(v)dv = \int_0^t du \, \theta^{\e,1}(u)\int_u^t dv  \theta^{\e,2}(v)
(v+\e-u)^{H-\frac12},
$$
and the same kind of argument also gives
\begin{eqnarray*}
&&\int_0^t Y^{1,\e}(v,t)\theta^{\e,2}(v)dv\\
 &=&
\int_0^t du \,\theta^{\e,1}(u)\int_u^t dw\int_u^w dv\,\theta^{\e,2}(v)\,(w-v)^{H-\frac12}
\big( 
(w+\e-u)^{H-\frac12}-(v+\e-u)^{H-\frac12}
\big)
\\
&&+\int_0^t du\,\theta^{\e,1}(u)\int_u^t dw\int_0^u dv\,\theta^{\e,2}(v)\,(w-v)^{H-\frac32}
(w+\e-u)^{H-\frac12}.
\end{eqnarray*}
Therefore, integrating first with respect to the randomness contained in $\theta^{\e,1}$, one is allowed to write $\delta_\e = E(
\Phi_\e(Z^\e)\,e^{iu_2\int_0^t \theta^{\e,2}(v)dv}
)$
where, for $f\in L^1([0,t])$, we set
$$
\Phi_\e(f):=E\left( e^{i\int_0^t f(u)\theta^{\e,1}(u)du} \right),
$$
and where the process $Z^\e$ is defined by:
\begin{eqnarray}\label{eq:36}
Z^\e_u&:=& u_1 + u_3 X^{2,\e}_u + u_4 Y^{2,\e}(u,t) + u_5\int_u^t (v+\e-u)^{H-\frac12}\theta^{\e,2}(v)dv\notag\\
&&+u_6\int_u^t dw \int_u^w dv \theta^{\e,2}(v)
\,(w-v)^{H-\frac32}\big( 
(w+\e-u)^{H-\frac12}-(v+\e-u)^{H-\frac12}
\big)\notag\\
&&+u_6\int_u^t dw \int_0^u dv\, \theta^{\e,2}(v)
(w-v)^{H-\frac32} 
(w+\e-u)^{H-\frac12}.\label{zeps}
\end{eqnarray}
Hence setting now, for $f\in L^{2}([0,t])$,
$$
\Phi(f):=E\left( e^{i\int_0^t f(u)dW^1_u}\right) = {\rm exp}\left(-\frac12\int_0^t f^{2}(u)du\right),
$$
we have obtained the decomposition
$$
\delta_\e = E\left( \Phi(Z)e^{iu_2 W^2_t}\right)  + v_\e^a + v_\e^b
$$
where the process $Z$ is given by
\begin{eqnarray*}
Z_u&=&u_1 + u_3 B^2_u + u_4 Y^{2}(u,t) + u_5 \int_u^t (v-u)^{H-\frac12}dW^2_v\\
&&+u_6\int_u^t dw \int_u^w dW^2_v\,(w-v)^{H-\frac32}\big( 
(w-u)^{H-\frac12}-(v-u)^{\frac12}
\big)\\
&&+u_6\int_u^t dw \int_0^u dW^2_v\,(w-v)^{H-\frac32} 
(w-u)^{H-\frac12},\quad u\in[0,t],
\end{eqnarray*}
and with two remainders $v_\e^a,v_\e^b$ defined as:
\begin{eqnarray*}
v_\e^a &:=& E\left( \Phi_\e(Z^\e)e^{iu_2 \int_0^t \theta^{\e,2}(u)du}\right)  - 
E\left( \Phi(Z^\e)e^{iu_2 \int_0^t \theta^{\e,2}(u)du}\right)\\
v_\e^b &:=& E\left( \Phi(Z^\e)e^{iu_2 \int_0^t \theta^{\e,2}(u)du}\right) 
- E\left( \Phi(Z)e^{iu_2 W^2_t}\right).
\end{eqnarray*}

\smallskip

The convergence of $v_\e^b$ above is easily established:
using again the same strategy than in \cite[Proof of Theorem 1]{DJ} 
(namely reducing the problem to a convergence of Kac-Stroock's process to white noise itself via an approximation of Liouville's kernel by step functions), one has that
$$
\left(Z^\e,\int_0^t \theta^{\e,2}(u)du\right)
\,\,\,\underset{\e\to 0}{\overset{\rm Law}{\longrightarrow}}\,\,\,(Z,W^2_t).
$$
Note that the convergence in law in the last equation holds in the space $\mathscr{C}\times\R$, where
$\mathscr{C}=\mathscr{C}([0,t])$ denotes the space of continuous function
endowed with the uniform norm $\|\cdot\|_\infty$. In particular, it is readily checked that 
$\lim_{\e\to 0} v_\e^b=0$.

\smallskip

Now, it remains to prove that
$
\lim_{\e \to 0} v_\e^a = 0.$ To this aim, we notice that we can bound trivially $|e^{iu_2 W^2_t}|$ by 1, and then apply the forthcoming Lemma \ref{lm-control2} in order to deduce that
$$
\vert v_{\e}^a  \vert \leq E  \left[ \left(\e^{2\alpha}
c_\alpha\,\|Z_{\e}\|_\alpha\|Z_{\e}\|_{L^2} u^2 +
\phi_{Z_{\e}}(\e) \frac{u^2}{2} + \psi_{Z_{\e}}(\e) \frac{u^4}{8}
+ \varphi_{Z_{\e}}(\e)\,\frac{|u|}{2}\right)e^{\frac{u^2\|Z_{\e}\|^2_{L^2}}{2}}\right]
$$
for any $\alpha\in(0,1)$. Furthermore, it is well known that characteristic functions on a neighborhood of 0 are sufficient to identify probability laws.
Consequently, using H\"older's inequality, we see that in order to get $\lim_{\e \to 0} v_\e^a = 0$, we are left to check that, for a given $u_0>0$,
\begin{eqnarray}
&&\sup_{0<\e\leq 1} E \big[ \|Z_{\e}\|_\alpha^2 \big] < \infty,
\label{next1}\\
&&\lim_{\e \to 0} E \big[ \phi^2_{Z_{\e}}(\e) \big] =0,
\quad
\lim_{\e \to 0} E \big[ \psi^2_{Z_{\e}}(\e) \big] =0,
\quad
\lim_{\e \to 0} E \big[ \varphi^2_{Z_{\e}}(\e) \big] =0,
\label{next4}\\
&&\sup_{0<\e\leq 1} E \big[  e^{{u^2\|Z_{\e}\|^2_{L^2}}}\big] \leq
M \quad \mbox{ for all $u \leq u_0$}\label{next5}.
\end{eqnarray}
We are now going to see that relations (\ref{next1}), (\ref{next4}) and (\ref{next5}) are satisfied.

\smallskip

\noindent
{\it (iii) Simplification of inequality (\ref{next5})}. Recall that $Z^\e$ has been defined by (\ref{eq:36}), and decompose it as $Z^{\e}_u=u_1 + u_3 U_1^{\e} (u) + u_4 U_2^{\e} (u) +u_5 U_3^{\e}
(u) + u_6 U_4^{\e} (u) + u_6 U_5^{\e} (u)$, with 
\begin{align*}
& U_1^{\e} (u)=X_u^{2,\e}, \quad U_2^{\e} (u)=Y^{2,\e} (u,t),
\quad U_3^{\e}(u)=\int_u^t (r+\e-u)^{H-\frac12} \theta^{\e,2} (r) dr \\
& U_4^{\e} (u)=u_6 \int_u^t dw \int_u^w dr \, \theta^{\e,2}(r)
 (w-r)^{H-\frac32} 
( (w+\e-u)^{H-\frac12} - (r+\e-u)^{H-\frac12})  \\
&U_5^{\e} (u)=\int_u^t dw \int_0^u dr \, \theta^{\e,2}(r)
(w-r)^{H-\frac32}
(w+\e-u)^{H-\frac12}
\end{align*}
In order to obtain (\ref{next5}), it is sufficient to check that there exists $M>0$ such that, for $\kappa>0$ small enough and 
$i=1,\ldots,5$, we have
\begin{equation}\label{etoile}
\sup_{0<\e\leq 1} E \left(  e^{{\kappa \int_0^T U_i^{\e} (u)^2
du }}\right) \leq M.
\end{equation}
Moreover, observe that $U_i^\e$ can always be written under the form
\begin{equation}\label{vi}
U_i^{\e} (u) = \int_0^T V_i(u,r,\e) \theta^{\e,2} (r) dr,
\end{equation} 
for a deterministic function $V_i(u,r,\e)$, and it is thus enough to check that
\begin{equation}\label{eq:42}
C_i:=\sup_{u\in[0,T]}\,\,\sup_{0<\e\leq 1} \int_0^T V_i^2(u,r,\e) dr < \infty.
\end{equation}
Indeed, using Lemma \ref{lm-control} below, we can write
\begin{eqnarray*}
 & & E \left(  e^{{\kappa \int_0^T U_i^{\e} (u)^2
du }}\right) = \sum_{m=0}^\infty \frac{\kappa^m}{m!} E \left[
\left(
 \int_0^T U_i^{\e} (u)^2 du \right)^m \right]
  \\
&   \leq&\frac1T \sum_{m=0}^\infty \frac{(T\kappa)^m}{m!}  \int_0^T E
\big[U_i^{\e} (u)^{2m}\big] du\\
& \leq&\frac1T \int_0^T \sum_{m=0}^\infty \frac{(2m)!(T\kappa)^m}{2^m (m!)^2}
 \Vert V_i(u,\cdot,\e) \Vert_{L^2}^{2m} du 
  \leq \sum_{m=0}^\infty ( 9 T\kappa C_i)^m,
\end{eqnarray*}
where we have used the bound $(m/3)^m\le m!\le m^m$ in the last inequality,
so that the desired conclusion follows for $\kappa>0$ small enough.

\smallskip

\noindent
{\it (iv) Proof of (\ref{eq:42}).}
We shall treat separately the cases for $i=1,\ldots,5$. During all the computations below, $C>0$ will denote a constant
depending only on $H$ and $T$, which can differ from one line to another.

\smallskip

\noindent
a) {\it Case $i=1$}. We have
$
X_u^{2,\e}  = \int_0^T
V_1(u,r,\e) \theta^{\e,2} (r) dr
$
with
$$
V_1(u,r,\e) = {\bf 1}_{[0,u]}(r) (u+\e-r)^{H-\frac12}
$$
Since
$$
\int_0^T V_1^2(u,r,\e) dr = \int_0^u (u+\e-r)^{2H-1} dr \leq 
\int_0^{u} (u-r)^{2H-1} dr = \frac{u^{2H}}{2H}\leq C,$$
we have that (\ref{eq:42}) takes place for $i=1$.

\smallskip

\noindent
b)
{\it Case $i=2$}. We have
$
Y^{2,\e} (u,t) = \int_0^T V_2(u,r,\e) \theta^{\e,2}
(r) dr,
$
with
\begin{eqnarray*}
 V_2(u,r,\e) & = & {\bf 1}_{[0,u]}(r) \int_u^t \big( (w+\e-r)^{H-\frac12} - (u+\e-r)^{H-\frac12}\big) 
(w-u)^{H-\frac32} dw \\
 & + & {\bf 1}_{[u,t]}(r) \int_r^t  (w+\e-r)^{H-\frac12} (w-u)^{H-\frac32} dw .
\end{eqnarray*}
Then
$
\int_0^T V_2^2(u,r,\e) dr =A_{2,1}(u,\e) + A_{2,2}(u,\e),
$
where
\begin{eqnarray*}
A_{2,1} (u,\e)& = &  \int_0^u \left( \int_u^t \big( (w+\e-r)^{H-\frac12} -
(u+\e-r)^{H-\frac12}\big)
(w-u)^{H-\frac32} dw \right)^2 dr,\\
A_{2,2} (u,\e)& = &  \int_u^t \left( \int_r^t (w+\e-r)^{H-\frac12}
(w-u)^{H-\frac32} dw \right)^2 dr.
\end{eqnarray*}
For any $\beta \in (0,1)$ and $w>u>r>0$, we can write, for some $w^* \in (u+\e,w+\e)$:
\begin{eqnarray*}
& &\big\vert (w+\e-r)^{H-\frac12}- (u+\e-r)^{H-\frac12} \big\vert \\
& & \qquad   \leq \frac{C\vert w - u \vert^\beta}{\vert w^* - r
\vert^{(\frac32 -H)\beta}}
\left( \frac{1}{\vert w+\e-r \vert^{\frac12 -H}} + \frac{1}{\vert u+\e-r \vert^{\frac12 -H}}  \right)^{1-\beta}
  \leq  \frac{C\vert w - u \vert^\beta}{\vert u - r
\vert^{\frac12 +\beta-H}}.
\end{eqnarray*}
Then, choosing $\beta=\frac12 - H +
\delta$ (with $\delta>0$ small enough), we can write
\begin{eqnarray*}
A_{2,1}(u,\e)
& \leq &  C  \int_0^u    \frac{dr}{\vert u - r \vert^{2-4H+2\delta}}\times
\left( \int_u^t \frac{dw}{\vert w - u \vert^{1 -\delta}}
 \right)^2 \leq C,
\end{eqnarray*}
where we have used the fact that $2-4H<1$ whenever $H>1/4$.
On the other hand, using that, for $\alpha \in (0,1), \nu>1$ and
$\lambda>0$, we have
\begin{eqnarray*}
\int_0^\infty \frac{dy}{y^\alpha (y + \lambda)^\nu} 
= \frac{1}{\lambda^{\alpha+\nu-1}}
 \int_0^{\infty}
\frac{dy}{y^\alpha (y + 1)^\nu}
=
\frac{C}{\lambda^{\alpha+\nu-1}},
\end{eqnarray*}
we obtain that
\begin{eqnarray*}
A_{2,2}(u,\e) & \leq &   C \int_u^t \left( \int_r^t  \frac{dw}{\vert w-r \vert^{\frac12 -H}
\vert w-u \vert^{\frac32 -H}}  \right)^2 dr\\
&\leq &   C \int_u^t \left( \int_0^{\infty}  \frac{dy}{y^{\frac12 -H}(y+r-u)^{\frac32 -H}} \right)^2 dr \leq   
C \int_u^t   \frac{dr}{(r-u)^{2 -4H}} \leq C.
\end{eqnarray*}

\smallskip

\noindent
c) {\it Case $i=3$}. We have
\begin{equation*}
\int_u^t (r+\e-u)^{H-\frac12} \theta^{\e,2} (r) dr = \int_0^T V_3(u,r,\e)
\theta^{\e,2} (r) dr,
\end{equation*}
with
$
V_3(u,r,\e) = {\bf 1}_{[u,t]}(r) (r+\e-u)^{H-\frac12},
$
so that the desired conclusion follows immediately since
$$
\int_0^1 V_3^2(u,r,\e) dr = \int_u^t (r+\e-u)^{2H-1}  dr \leq \int_u^t (r-u)^{2H-1}dr =\frac{(t-u)^{2H}}{2H}\leq C.
$$

\smallskip

\noindent
d) {\it Case $i=4$}. We can write
\begin{equation*}
\int_u^t dw \int_u^w dr (w-r)^{H-\frac32} \big(
(w+\e-u)^{H-\frac12} - (r+\e-u)^{H-\frac12}\big)\theta^{\e,2} (r)  
\end{equation*}
as $\int_0^T V_4(u,r,\e) \theta^{\e,2} (r) dr$, with
\begin{equation*}
V_4(u,r,\e) = {\bf 1}_{[u,t]}(r) \int_r^t  
(w-r)^{H-\frac32} \big( (w+\e-u)^{H-\frac12} - (r+\e-u)^{H-\frac12}\big) dw.
\end{equation*}
Then, according to the computations already performed for the analysis of 
$A_{2,1}$ above, we obtain, for $\delta>0$ small enough,
\begin{equation*}
\int_0^T V_4^2(u,r,\e) dr   \leq   C  \int_u^t \frac{1}{\vert r -
u \vert^{2-4H+2\delta}} \left( \int_r^t \frac{dw}{\vert w - r
\vert^{1 -\delta}} \right)^2 dr \leq C.
\end{equation*}

\smallskip

\noindent
e) {\it Case $i=5$}. We have
\begin{equation*}
\int_u^t dw \int_0^u dr  (w-r)^{H-\frac32}
(w+\e-u)^{H-\frac12} \theta^{\e,2} (r)  = \int_0^T V_5(u,r,\e)
\theta^{\e,2} (r)dr,
\end{equation*}
with
$$
V_5(u,r,\e) = {\bf 1}_{[0,u]}(r) \int_u^t  (w-r)^{H-\frac32}  (w+\e-u)^{H-\frac12} dw.
$$
Since $\vert w- r \vert^{\frac32-H} \geq \vert w - u
\vert^{1-H+\delta} \vert u - r \vert^{\frac12 - \delta}$ for $r<u<w$, we get
(for $\delta>0$ small enough) that
\begin{eqnarray*}
\left\vert (w-r)^{H-\frac32}  (w+\e-u)^{H-\frac12} \right\vert &
\leq & \frac{C}{\vert w - r \vert^{\frac32 - H}\vert w - u \vert^{\frac12 - H}} 
 \leq  \frac{C}{\vert u - r \vert^{\frac12 - \delta}\vert w - u \vert^{\frac32 - 2H+\delta}}.
\end{eqnarray*}
Hence, invoking again the fact that $H>1/4$, we end up with
\begin{equation*}
\int_0^T V_5^2(u,r,\e) dr   \leq   C  \int_0^u \frac{dr}{\vert u -
r \vert^{1 - 2\delta}}\times \left( \int_u^t \frac{dw}{\vert w - u
\vert^{\frac32 - 2H+\delta}}
  \right)^2  \leq C.
\end{equation*}

\smallskip

\noindent
{\it (v) Proof of (\ref{next4}).} In the previous step, we have
shown in particular that, for any $i=1,...,5$, we have 
$\sup_{0<\e\leq 1} \int_0^T E\big[ |U_i^{\e} (u)|^p \big]du  <\infty$ for all $p>1$,
which implies
$$
\sup_{0<\e\leq 1} \int_0^T E\big[ |Z^{\e}_u|^p \big]du  <\infty,\quad\mbox{for all $p>1$}.
$$
On the other hand, a simple application of Schwarz inequality yields
\begin{equation*}
 E \big[ \phi^2_{Z^{\e}}(\e) \big] =
 E \left[ \left( \int_0^T (Z^{\e}_u)^2  e^{-\frac{2u}{\e^2}} du \right)^2\right]
 \leq C \e^2  \int_0^T E \big[ (Z^{\e}_u)^4 \big] du,
 \end{equation*}
 and the same kind of argument also gives
\begin{align*}  &E \big[ \psi^2_{Z^{\e}}(\e) \big]  =
 E \left[ \left( \int_0^T dx \int_0^x  dy (Z^{\e}_x)^2 (Z^{\e}_y)^2
e^{-\frac{2(x-y)}{\e^2}} \right)^2\right]\\
& \leq  \frac12 E \left[ \left( \int_0^T dx \int_0^x  dy
(Z^{\e}_x)^4 e^{-\frac{2(x-y)}{\e^2}} \right)^2\right] +
\frac12 E \left[ \left( \int_0^T dx \int_0^x  dy  (Z^{\e}_y)^4
e^{-\frac{2(x-y)}{\e^2}} \right)^2\right]\\
& \leq  C \e^4   \int_0^T E \big[ (Z^{\e}_u)^8 \big] du.
\end{align*}
Finally, we have
\begin{align*}  
&E \big[ \varphi^2_{Z^{\e}}(\e) \big]  =
E \left[  \left( \e \Vert Z^{\e} \Vert_{L^2} + \lp \int_0^\e
(Z^{\e}_u)^2 du \rp^{1/2} \right)^2 \right]\\
& \leq  2  \e^2 E ( \Vert Z^{\e} \Vert_{L^2}^2 ) +
2 \int_0^\e E\big[ ( Z^{\e}_u )^2 \big] du \leq  2  \e^2 E ( \Vert Z^{\e} \Vert_{L^2}^2 ) + 2
\e^{1/2}\, \lp \int_0^T E\big[ (Z^{\e}_u)^4 \big] du\rp^{1/2},
\end{align*}
and the proof of (\ref{next4}) follows immediately by putting all these facts together.\\

\smallskip

\noindent
{\it (vi) Proof of (\ref{next1})}.
For all $\alpha < \beta
- \frac{1}{p}$, the Sobolev inequality (\ref{eq:24}) yields $\Vert Z^{\e} \Vert_\alpha \leq
C\Vert  Z^{\e} \Vert_{\beta,p}$, where $\|f\|_{\beta,p}$ has been defined by (\ref{eq:def-sobolev-norm}).
Moreover, recall from (\ref{zeps}) that $Z^\e$ has the form
$$
Z^\e_t -Z^\e_s = \int_0^T G(s,t,r) \theta^{\e,2}(r)dr
$$
for some $G(s,t,\cdot)\in L^2([0,T])$.
Hence, using the definition of $\theta^{\e,2}$, we can write, for any even integer $p\geq 2$,
\begin{eqnarray*}
&&E|Z^\e_t-Z^\e_s|^p\\
&=&\e^{-p}\int_{[0,T]^p}G(s,t,r_1)\cdots G(s,t,r_p) E\big[(-1)^{N(\frac{r_1}{\e})
+\ldots+N(\frac{r_p}{\e})}\big] dr_1\cdots dr_p\\
&=&p!\e^{-p}\int_{[0,T]^p}G(s,t,r_1)\cdots G(s,t,r_p) 
e^{-\frac{2(r_1-r_2)}{\e^2}}
\cdots
e^{-\frac{2(r_{p-1}-r_p)}{\e^2}}{\bf 1}_{\{r_1\geq\cdots\geq r_p\}}
 dr_1\cdots dr_p\\
&=&
\frac{p!}{(p/2)!}
\left(\e^{-2}\int_{[0,T]^2}G(s,t,r_1)G(s,t,r_2) 
e^{-\frac{2(r_1-r_2)}{\e^2}}
{\bf 1}_{\{r_1\geq r_2\}}
 dr_1dr_2\right)^{p/2}\\
&=&\frac{p!}{(p/2)!}
\left(\frac{\e^{-2}}{2}\int_{[0,T]^2}G(s,t,r_1)G(s,t,r_2) 
E\big[(-1)^{N(\frac{r_1}{\e})+N(\frac{r_2}{\e})}\big] 
 dr_1dr_2\right)^{p/2}\\
&=&\frac{p!}{2^{\frac{p}2}(p/2)!}\big(E|Z^\e_t-Z^\e_s|^2\big)^{p/2}.
\end{eqnarray*}
In particular, we see that, in order to achieve the proof of (\ref{next1}), it is enough to check that
\begin{equation}\label{needed}
E \vert Z^{\e}_t - Z^{\e}_s \vert^2 \leq C \vert t - s \vert^{H-\delta}
\end{equation}
for some $\delta>0$ (small enough).
Actually, we shall use again the decomposition of $Z^\e$ in terms of the $U_i$'s, which means that it is sufficient to prove $E \vert U_i^{\e}(u) - U_i^{\e}(v) \vert^2 \leq C \vert u-v \vert^{H-\delta}$ for $i=1,...,5$. But it is easily seen that
\begin{eqnarray*}
E \vert U_i^{\e} (u)- U_i^{\e} (v) \vert^2 &=&
E \left|
 \int_0^T \big(V_i(u,r,\e) -  V_i(v,r,\e) \big)\theta^{\e,2} (r) dr \right|^2 \\
 & \leq&    \int_0^T \big(V_i(u,r,\e) -  V_i(v,r,\e) \big)^2 dr ,
\end{eqnarray*}
where $V_i$ is defined by (\ref{vi}), and are specified at step (iii). It is thus enough for our purposes to show that 
$ \int_0^T(V_i(u,r,\e) -  V_i(v,r,\e))^2 dr \leq C \vert u - v \vert^{H-
 \delta}$ for $0\le u\le v\le t$, where $t\in\ott$, and we will now focus on this last inequality, for $i=1,...,5$.

\smallskip

\noindent
a) {\it Case $i=1$}. We have
\begin{eqnarray*}
&& \int_0^T \big(V_1(u,r,\e) -  V_1(v,r,\e) \big)^2 dr\\
& & \quad = \int_0^v \big( (u+\e-r)^{H-\frac12} - (v+\e-r)^{H-\frac12} \big)^2 dr + \int_v^u
(u+\e-r)^{2H-1} dr \\
& & \quad \leq \int_0^v \big( (u-r)^{H-\frac12} - (v-r)^{H-\frac12} \big)^2 dr + \int_v^u
(u-r)^{2H-1} dr \\
& & \quad \leq c(u)^{2H}\quad\mbox{by Lemma \ref{onenabesoin}}.
\end{eqnarray*}

\smallskip

\noindent
b) {\it Case $i=2$}. We can decompose $\int_0^T \big(V_2(u,r,\e) -  V_2(v,r,\e) \big)^2 dr$ into
$$
 \int_0^t \big(V_2(u,r,\e) -  V_2(v,r,\e) \big)^2 dr=
 B_{2,1} + B_{2,2}+B_{2,3}+B_{2,4},
$$
with
\begin{eqnarray*}
B_{2,1}&=&\int_0^v \left( \int_u^t \big( (w+\e-r)^{H-\frac12} - (u+\e-r)^{H-\frac12} \big)
(w-u)^{H-\frac32} dw \right. \\ && \left.\qquad \quad  -
\int_v^t \big( (w+\e-r)^{H-\frac12} - (v+\e-r)^{H-\frac12}\big) 
(w-v)^{H-\frac32} dw   \right)^2 dr  \\
B_{2,2}&=&\int_v^u \left( \int_u^t \big( (w+\e-r)^{H-\frac12} - (u+\e-r)^{H-\frac12}\big)
(w-u)^{H-\frac32} dw   \right)^2 dr
\end{eqnarray*}
and
\begin{eqnarray*}
B_{2,3}&=&
\int_v^u \left( \int_r^t  (w+\e-r)^{H-\frac12}   (w-v)^{H-\frac32}  dw  \right)^2 dr \\
B_{2,4}&=&
\int_u^t \left( \int_r^t  (w+\e-r)^{H-\frac12}  
\big((w-u)^{H-\frac32} -  (w-v)^{H-\frac32} \big) dw  \right)^2 dr.
\end{eqnarray*}

\smallskip

We will now study those terms separately, starting from $B_{2,1}$: this first term can be further decomposed as $B_{2,1}\le c( B_{2,1,1} + B_{2,1,2}+B_{2,1,3})$, with
\begin{align*}
&B_{2,1,1}  \\
&=\int_0^v \Bigg( \int_u^t \big( (w+\e-r)^{H-\frac12} - (u+\e-r)^{H-\frac12}\big) 
\times\big(  (w-u)^{H-\frac32}
-  (w-v)^{H-\frac32}\big) dw   \Bigg)^2 dr
\end{align*}
and
\begin{eqnarray*}
B_{2,1,2}&=&  
\int_0^v \left( \int_u^t \big(  (u+\e-r)^{H-\frac12} - (v+\e-r)^{H-\frac12}\big)
(w-v)^{H-\frac32}   dw   \right)^2 dr \\
B_{2,1,3}&=&  
\int_0^v \left( \int_v^u  \big( (w+\e-r)^{H-\frac12} - (v+\e-r)^{H-\frac12}\big)
(w-v)^{H-\frac32} dw   \right)^2 dr. 
\end{eqnarray*}
Moreover, for $\delta>0$ small enough, and by choosing $\beta_1=H-\delta$ and $\beta_2=\frac{H}{2}-2\delta$, we can write
\begin{eqnarray*}
B_{2,1,1} & \leq  & C \int_0^v  \left( \int_u^t   \frac{1}{\vert v-r \vert^{\beta_1(\frac32 -H)}} \vert w - u \vert^{\beta_1} \frac{1}{\vert v-r \vert^{(\frac12 - H)(1-\beta_1)}}\right.
\\ &&\left. \hspace{4cm}  \times \frac{1}{\vert w-u \vert^{\beta_2(\frac52 -H)}} \vert u- v \vert^{\beta_2} 
\frac{1}{\vert w-u \vert^{(\frac32 - H)(1-\beta_2)}}  dw   \right)^2 dr \\
& \leq  & C \vert u - v \vert^{2\beta_2} \int_0^v   \frac{1}{\vert v-r \vert^{1+2 \beta_1-2H}}     \left( \int_u^t    \frac{1}{\vert w-u \vert^{\frac32 -H +\beta_2- \beta_1}}  dw \right)^2   dr\\
& =  & C \vert u - v \vert^{H-4 \delta} \int_0^v   \frac{1}{\vert v-r \vert^{1 -2\delta}}     \left( \int_u^t    \frac{1}{\vert w-u \vert^{\frac32 - \frac32 H - \delta}}  dw \right)^2   dr
 \leq   C \vert u - v \vert^{H-4 \delta},
\end{eqnarray*}
since $\frac32 - \frac32 H - \delta < 1$ whenever $H>1/3$. The term $B_{2,1,2}$ can be bounded in the same way:
for  $\delta>0$ small enough, and by choosing again $\beta=H-\delta$, we can write
\begin{eqnarray*}
B_{2,1,2} & \leq  & C \int_0^v  \left( \int_u^t   \frac{1}{\vert v-r \vert^{\beta(\frac32 -H)}} \vert u - v \vert^\beta \frac{1}{\vert v-r \vert^{(\frac12 - H)(1-\beta)}} \frac{1}{\vert w-v \vert^{\frac32 -H}} dw \right)^2   dr\\
& \leq  & C \vert u - v \vert^{H-4 \delta} \int_0^v   \frac{1}{\vert v-r \vert^{1-2\delta}}     \left( \int_v^t    \frac{1}{\vert w-v \vert^{\frac32 - \frac32 H - \delta}}  dw \right)^2   dr
 \leq   C \vert u - v \vert^{H-4 \delta}.
\end{eqnarray*}
Finally, for  $\delta>0$ small enough, and still with $\beta=H-\delta$, we can write
\begin{eqnarray*}
B_{2,1,3} & \leq  & C \int_0^v  \left( \int_v^u   \frac{1}{\vert v-r \vert^{\beta(\frac32 -H)}} \vert w - v \vert^\beta \frac{1}{\vert v-r \vert^{(\frac12 - H)(1-\beta)}} \frac{1}{\vert w-v \vert^{\frac32 -H}} dw \right)^2   dr\\
& \leq  & C \int_0^v   \frac{1}{\vert v-r \vert^{1-2\delta}}     \left( \int_v^u    \frac{1}{\vert w-v \vert^{\frac32 - 2H+\delta}}  dw \right)^2   dr\\
& \leq  & C \vert u - v \vert^{ 4H-1-2\delta} \int_0^v   \frac{1}{\vert v-r \vert^{1-2\delta}}      dr
\leq C|u-v|^{H-2\delta}.
\end{eqnarray*}
By gathering the bounds we have obtained for $B_{2,1,1},B_{2,1,2},B_{2,1,3}$, we now easily get $B_{2,1}\le C |u-v|^{H-\delta}$ for any arbitrarily small constant $\delta>0$.

\smallskip

Let us turn now to the bound concerning $B_{2,2}$: for  $\beta=\frac12 - H +\delta$ with $\delta>0$ small enough, we can write,
by using similar arguments to those used in the study of $B_{2,1,3}$,
\begin{eqnarray*}
B_{2,2} & \leq  & C \int_v^u  \left( \int_u^t  \frac{1}{\vert u-r
\vert^{\beta( \frac32 -H)}} \vert w - u \vert^\beta \frac{1}{\vert
u-r \vert^{(\frac12-H)(1-\beta)}}  \frac{1}{\vert w-u
\vert^{ \frac32 -H}} dw \right)^2 dr \\
& \leq & C    \int_v^u \frac{1}{\vert u-r \vert^{ 2-4H +2 \delta }}
\left( \int_u^t \frac{1}{\vert w-u \vert^{ 1 - \delta }}  dw
\right)^2 dr
 \leq  C   \vert u-v \vert^{H-2 \delta}.
\end{eqnarray*}
In a similar manner we obtain, still for $\delta>0$ small enough, 
\begin{eqnarray*}
B_{2,3} & \leq  & C \int_v^u \left(  \int_r^t  \frac{1}{\vert w-r
\vert^{\frac12 -H}} \frac{1}{\vert w-v \vert^{\frac32 -H}}
  dw  \right)^2  dr \\
& \leq &  C \int_v^u  \frac{1}{\vert v-r \vert^{2-4H+2 \delta}}
\left(  \int_r^t  \frac{1}{\vert w-r \vert^{1-  \delta}}
  dw  \right)^2  dr
 \leq  C    \vert u-v \vert^{H-2\delta}  .
\end{eqnarray*}
Finally, by choosing $\beta=\frac{H}{2} - 2\delta$ with $\delta>0$ small enough,
the same kind of arguments allow to write
\begin{eqnarray*}
B_{2,4} & \leq  & C \int_u^t  \left( \int_r^t   \frac{1}{\vert w-r
\vert^{ \frac12 -H}}  \frac{1}{\vert w-u \vert^{\beta( \frac52
-H)}} \vert u - v \vert^\beta \frac{1}{\vert
w-u \vert^{(\frac32-H)(1-\beta)}}  dw \right)^2 dr \\
& \leq & C   \vert u-v \vert^{H-4 \delta}  \int_u^t \frac{1}{\vert
r-u \vert^{ 2-3H -2 \delta }} \left( \int_r^t \frac{1}{\vert w-r
\vert^{ 1 - \delta }}  dw \right)^2 dr
 \leq  C   \vert u-v \vert^{H-4 \delta}.
\end{eqnarray*}
In conclusion, putting together the bounds we have obtained for all the terms $B$, we end up with the announced claim $\int_0^T \big(V_2(u,r,\e) -  V_2(v,r,\e) \big)^2 dr\le C |u-v|^{H-\delta}$ for any  $0<\delta<H$, as soon as $H>1/3$.

\smallskip

\noindent
c) {\it Case $i=3$}. This case can be treated along the same lines as the previous ones: we have
\begin{eqnarray*}
&& \int_0^T \big(V_3(u,r,\e) -  V_3(v,r,\e) \big)^2 dr\\
& & \quad = \int_0^T \big({\bf 1}_{[u,t]}(r) (r+\e-u)^{H-\frac12} - {\bf 1}_{[v,t]}(r) (r+\e-v)^{H-\frac12} \big)^2  dr \\
& & \quad \leq C \left(  \int_v^u (r-v)^{2H-1} dr + \int_u^t  \big(
(r-u)^{H-\frac12} - (r-v)^{H-\frac12}\big)^2 dr
\right).
\end{eqnarray*}
Now it can be shown, by using the same reasoning as in the proof of Lemma \ref{onenabesoin}, that $\int_0^1 \big(V_3(u,r,\e) -  V_3(v,r,\e) \big)^2 dr\le C|u-v|^{2H}$.

\smallskip

\noindent
d) {\it Case  $i=4$}. Let us decompose again the integral under consideration into:
$$
\int_0^T \big(V_4(u,r,\e) -  V_4(v,r,\e) \big)^2 dr=
B_{4,1} + B_{4,2},
$$
with
\begin{eqnarray*}
B_{4,1}&=&
\int_v^u \left(  \int_r^t  (w-r)^{H-\frac32} \big( (w+\e-v)^{H-\frac12} - 
(r+\e-v)^{H-\frac12}\big) dw  \right)^2  dr  \\
B_{4,2}&=&
\int_u^t \left(  \int_r^t  (w-r)^{H-\frac32} \big( 
(w+\e-u)^{H-\frac12} - (r+\e-u)^{H-\frac12} \right.  \\
& & \hskip6cm \left.- (w+\e-v)^{H-\frac12} + (r+\e-v)^{H-\frac12}\big) 
dw  \right)^2  dr .
\end{eqnarray*}
The first of these terms can now be treated by using similar arguments to those used in the study of $B_{3,2}$: for  $\beta=\frac12 - H +\delta$ (with $\delta>0$ small enough), we can write
\begin{eqnarray*}
B_{4,1} & \leq  & C \int_v^u  \left( \int_r^t \frac{1}{\vert w-r
\vert^{ \frac32 -H}}  \frac{1}{\vert r-v \vert^{\beta( \frac32
-H)}} \vert w - r \vert^\beta \frac{1}{\vert v-r
\vert^{(\frac12-H)(1-\beta)}} dw \right)^2 dr \\
& \leq & C    \int_v^u \frac{1}{\vert v-r \vert^{ 2-4H +2 \delta }}
\left( \int_r^t \frac{1}{\vert w-r \vert^{ 1 - \delta }}  dw
\right)^2 dr
 \leq  C   \vert u-v \vert^{H-2 \delta}.
\end{eqnarray*}

\smallskip

In order to estimate $B_{4,2}$, we introduce an interpolated point 
$(x^*,y^*) \in (r+\e, w+\e) \times (v,u)$ which allows to write, for $u<r$,
\begin{eqnarray*}
&& \big| (w+\e-u)^{H-\frac12} - (r+\e-u)^{H-\frac12}- (w+\e-v)^{H-\frac12} + (r+\e-v)^{H-\frac12} \big| \\
& & \quad \leq \left(\frac12-H\right)\left(\frac32-H\right) |x^*-y^*|^{H-\frac52} \vert w - r \vert \vert u - v \vert \leq  C \,
\frac{\vert w - r \vert
\vert u - v \vert}{\vert r-u \vert^{ \frac52 -H}}  .
\end{eqnarray*}
Therefore, with $\beta=\frac{H}{2}-\delta$ and $\delta>0$ such that
$\frac32 -\frac32 H + \delta <1$ (recall again that this is possible as soon as $H >
\frac13$), we get
\begin{eqnarray*}
B_{4,2} & \leq  & C \int_u^t  \left( \int_r^t \frac{1}{\vert w-r
\vert^{ \frac32 -H}}  \frac{1}{\vert r-u \vert^{ \beta (\frac52
-H)}} \vert \vert w - r \vert^\beta  \vert u - v \vert^\beta
\frac{1}{\vert r-u
\vert^{(\frac12-H)(1-\beta)}} dw \right)^2 dr \\
& \leq & C    \vert u - v \vert^{2\beta} \int_u^t \frac{1}{\vert
r-u \vert^{1+4\beta - 2 H} } \left( \int_r^t \frac{1}{\vert w-r
\vert^{ \frac32 -H - \beta}}    dw \right)^2 dr\\
& \leq & C    \vert u - v \vert^{H-2\delta} \int_u^t \frac{1}{\vert
r-u \vert^{1 - 2 \delta} } \left( \int_r^t \frac{1}{\vert w-r
\vert^{ \frac32 -\frac32 H + \delta}}    dw \right)^2 dr
 \leq  C   \vert u-v \vert^{H-2 \delta}.\\
\end{eqnarray*}

\smallskip

\noindent
e) {\it Case $i=5$}. Our last decomposition can be written $\int_0^1 \big(V_5(u,r,\e) -  V_5(v,r,\e) \big)^2 dr\le C(B_{5,1} + B_{5,2}+B_{5,3})$, with 
\begin{eqnarray*}
B_{5,1}&=&  
\int_v^u \left(  \int_u^t (w-r)^{H-\frac32}  (w+\e-u)^{H-\frac12} dw  \right)^2  dr \\
B_{5,2}&=&  
\int_0^v \left(  \int_u^t (w-r)^{H-\frac32} \big( (w+\e-u)^{H-\frac12} - (w+\e-v)^{H-\frac12} 
\big) dw   \right)^2  dr \\
B_{5,3}&=&  
\int_0^v  \left(
\int_v^u  (w-r)^{H-\frac32}  (w+\e-v)^{H-\frac12} dw   \right)^2  dr.
\end{eqnarray*}

Then, choosing $\delta>0$ small enough, and thanks to the fact that 
\begin{equation}\label{eq:46}
(w-r)^{3/2-H}=(w-r)^{1-2H+\delta}(w-r)^{1/2+H-\delta}\geq (u-r)^{1-2H+\delta}(w-u)^{1/2+H-\delta}
\end{equation}
for all $v<r<u<w<t$, we have
\begin{multline*}
B_{5,1}  \leq   C
\int_v^u \left(  \int_u^t  \frac{dw}{\vert
w-r \vert^{\frac32 -H}\vert
w-u \vert^{\frac12 -H}}
   \right)^2  dr \\
 \leq   C
\int_v^u  \frac{1}{\vert
u-r \vert^{2-4H+2 \delta}}   \left(  \int_u^t  \frac{1}{\vert
w-u \vert^{1-\delta}}
  dw  \right)^2  dr
 \leq  C    \vert u-v \vert^{4H-1-2\delta}  \leq  C   \vert u-v \vert^{H-2\delta}  .
\end{multline*}
On the other hand, for  $\beta=\frac{H}{2}-2\delta$ with $\delta>0$ small enough, we can write,
using similar arguments to those used in the study of $B_{2,3}$,
\begin{eqnarray*}
B_{5,2} & \leq  & C \int_0^v  \left( \int_u^t \frac{1}{\vert
w-r \vert^{ \frac32 -H}}  \frac{1}{\vert
w-u \vert^{\beta( \frac32 -H)}} \vert u - v \vert^\beta \frac{1}{\vert w-u
\vert^{(\frac12-H)(1-\beta)}} dw \right)^2 dr \\
& \leq & C   \vert u-v \vert^{H-4 \delta}    \int_0^v   \frac{1}{\vert
v-r \vert^{ 1 -2 \delta }}   \left( \int_u^t  \frac{1}{\vert
w-u \vert^{ \frac32 -\frac32 H - \delta }}  dw \right)^2 dr
 \leq  C   \vert u-v \vert^{H-4 \delta}.
\end{eqnarray*}
Finally, still for $\delta>0$ small enough, invoking again inequality (\ref{eq:46}), we have 
\begin{eqnarray*}
B_{5,3} & \leq  & C
\int_0^v \left(  \int_v^u  \frac{1}{\vert
w-r \vert^{\frac32 -H}}
\frac{1}{\vert
w-v \vert^{\frac12 -H}}
  dw  \right)^2  dr \\
& \leq &  C
\int_0^v  \frac{1}{\vert
v-r \vert^{2-3H-2 \delta}}   \left(  \int_v^u  \frac{1}{\vert
w-v \vert^{1- \frac{H}{2} + \delta}}
  dw  \right)^2  dr
  \leq  C    \vert u-v \vert^{H-2\delta}  .
\end{eqnarray*}
Our estimates on $B_{5,1},B_{5,2},B_{5,3}$ yield trivially $\int_0^1 \big(V_5(u,r,\e) -  V_5(v,r,\e) \big)^2 dr\le C \vert u-v \vert^{H-\delta}$ for $\delta$ arbitrarily close to 0, and this last bound finishes the proof of our Proposition \ref{prop-weak}.

\end{proof}

\section{Some technical lemmas}\label{sec:tech}
This section collect the technical results that have been used throughout the proof of 
Theorem \ref{main-thm-weak}.
The first lemma aims at giving some 
estimates concerning the Kac-Stroock kernel (\ref{stroock}), which can be seen as a elaboration of  the ones contained in Delgado and Jolis \cite[Lemma 2]{DJ}. Notice however that these latter results are not sharp enough for our purposes, which forced us to a refinement.

\begin{lemma}\label{lm-control}
Let $m\in\N$, $f,f_1,\ldots,f_{2m}\in L^2([0,T])$,
$k\in\{1,2\}$ and $\e>0$. We have:
\begin{equation}\label{lm-cont-pair}
\left|E\left[\prod_{j=1}^{2m} \int_0^T f_j(r)\theta^{\e,k}(r)dr
\right]\right|\leq\frac{(2m)!}{2^mm!}\,\|f_1\|_{L^2}\ldots\|f_{2m}\|_{L^2},
\end{equation}
and
\begin{equation}\label{lm-cont-impair}
\left|E\left[\left( \int_0^T f(r)\theta^{\e,k}(r)dr \right)^{2m+1}
\right]\right|\leq\varphi_f(\e)\,\frac{(2m+1)!}{2^{m+1}m!}\,\|f\|_{L^2}^{2m},
\end{equation}
where
$$
\varphi_f(\e)=\e \|f\|_{L^2} + \left( \int_0^\e \vert f(s) \vert^2
ds \right)^\frac12.
$$
\end{lemma}
\begin{proof}
For $m\in\N$, $\e>0$ and $f_1,\ldots,f_{2m}\in L^2([0,T])$,
let us denote
$$
\Delta^\e_{2m}(f_1,\ldots,f_{2m}):=E\left[\prod_{j=1}^{2m}
\int_0^T f_j(r)\theta^{\e,k}(r)dr\right].
$$

\smallskip

We will need to introduce some operations on the set of permutations (in the sequel, $\mathfrak{S}_{k}$ stands for the set of permutations on $\{1,\ldots,k\}$):
when $\tau\in\mathfrak{S}_{2m}$ and
$\sigma\in\mathfrak{S}_m$, we note $\sigma\star\tau$ the element of $\mathfrak{S}_{2m}$
defined by
$$
(\sigma\star\tau)(2j-1)=\tau(2\sigma(j)-1)\quad\mbox{and}\quad(\sigma\star\tau)(2j)=\tau(2\sigma(j)).
$$
Remark that we have ${\rm id}\star\tau=\tau$ and $\sigma'\star(\sigma\star\tau)=(\sigma\sigma')\star\tau$, so
$\star:\mathfrak{S}_{m}\times\mathfrak{S}_{2m}\rightarrow\mathfrak{S}_{2m}$ defines a (right) group action of $\mathfrak{S}_m$ on $\mathfrak{S}_{2m}$.
For any $\tau\in\mathfrak{S}_{2m}$, the orbit of $\tau$ has exactly $m!$ elements.
Consequently, the set $\mathscr{O}$ of the orbits under the group action $\star$ has $\frac{(2m)!}{m!}$ elements and
we have, by denoting $\tau_i$ one particular element of the orbit $o_i=o(\tau_i)\in\mathscr{O}$:
for $r_1,\ldots,r_{2m}\in[0,1]$,
\begin{eqnarray}\label{eq:47}
{\bf 1}_{\{\forall i\neq j,\,\,r_i\neq r_j\}}
=\sum_{\tau\in\mathfrak{S}_{2m}} {\bf 1}_{\{r_{\tau(1)}>\ldots>r_{\tau(2m)}\}}
&=&\sum_{o_i\in\mathscr{O}} \sum_{\tau\in o_i}{\bf 1}_{\{r_{\tau(1)}>\ldots>r_{\tau(2m)}\}}
\notag\\
&\leq&\sum_{i=1}^{\nicefrac{(2m)!}{m!}} \prod_{j=1}^m {\bf 1}_{\{r_{2\tau_i(j)-1}>r_{2\tau_i(j)}\}}.
\end{eqnarray}
For the reader who might not be completely convinced by this inequality, let us illustrate it by an example: when $m=2$ and $\tau_i={\rm id}\in\mathfrak{S}_4$, we have $o_i=o(\tau_i)=\{{\rm id},(13)(24)\}$
and we have used
\begin{eqnarray*}
\sum_{\tau\in o_i}{\bf 1}_{\{r_{\tau(1)}>r_{\tau(2)}>r_{\tau(3)}>r_{\tau(4)}\}}
={\bf 1}_{\{r_1>r_2>r_3>r_4\}}+{\bf 1}_{\{r_3>r_4>r_1>r_2\}}&\leq& {\bf 1}_{\{r_1>r_2\}}{\bf 1}_{\{r_3>r_4\}}\\
&=&\prod_{j=1}^2 {\bf 1}_{\{r_{2\tau_i(j)-1}>r_{2\tau_i(j)}\}}.
\end{eqnarray*}

\smallskip

Let us apply now inequality (\ref{eq:47}). We introduce first a notation which will prevail until the end of the article: for $\e>0$ and $r\in\R_+$, we set $Q_{\e}(r):= e^{-\frac{2r}{\e^2}}/\e^2$.
Notice then that, for any $\e>0$:
\begin{align*}
&|\Delta^\e_{2m}(f_1,\ldots,f_{2m})|\\
&\leq \frac{1}{\e^{2m}}\int_{[0,T]^{2m}}
|f_1(r_1)|\ldots |f_{2m}(r_{2m})| \left|E\left[(-1)^{
\sum_{i=1}^{2m} N(\frac{r_{i}}{\e^2}) }\right]\right|
dr_1\ldots dr_{2m}\\
&=\sum_{o_i\in\mathscr{O}} \sum_{\tau\in o_i}
\frac{1}{\e^{2m}}
\int_{[0,T]^{2m}}
\!\!\!\!\!{\bf 1}_{\{r_{\tau(1)}>\ldots>r_{\tau(2m)}\}}
|f_1(r_1)|\cdots |f_{2m}(r_{2m})| \\
&\hspace{8cm} Q_{\e}\Bigg(  \sum_{i=1}^m (r_{\tau(2i-1)}-r_{\tau(2i)})\Bigg)
dr_1\cdots dr_{2m}
\end{align*}
and thus, according to (\ref{eq:47}), we obtain
\begin{align*}
&|\Delta^\e_{2m}(f_1,\ldots,f_{2m})| \\
&\leq\sum_{i=1}^{\nicefrac{(2m)!}{m!}}
\prod_{j=1}^m
\int_{[0,T]^2} {\bf 1}_{\{r_1> r_2\}} |f_{2\tau_i(j)-1}(r_1)||f_{2\tau_i(j)}(r_2)|
Q_\e(r_{1}-r_{2})
dr_1dr_2  \\
&\leq\frac{(2m)!}{2^mm!}\|f_1\|_{L^2}\ldots\|f_{2m}\|_{L^2},
\end{align*}
the last inequality coming from
\begin{eqnarray}
& &\int_{[0,T]^2} \!\!\!\!{\bf 1}_{\{r_1> r_2\}}
|f_{k}(r_1)||f_{\ell}(r_2)| \, Q_\e(r_{1}-r_{2}) dr_1dr_2 \nonumber\\
&\leq& \left( \int_0^T \!\!\!\! |f_{k}(r_1)|^2  \left( \int_0^{r_1}
Q_\e(r_{1}-r_{2}) dr_2
\right) dr_1 \right)^\frac12 \left( \int_0^T \!\!\!\!
|f_{\ell}(r_2)|^2  \left( \int_{r_2}^T Q_\e(r_{1}-r_{2}) dr_1 \right) dr_2
\right)^\frac12\nonumber
\\
&\leq&\frac{1}{2}\|f_k\|_{L^2}\|f_\ell\|_{L^2}.\label{ineone}
\end{eqnarray}

\smallskip

This finishes the proof of (\ref{lm-cont-pair}), so let us now concentrate on (\ref{lm-cont-impair}).
For $m\in\N$, we have
\begin{eqnarray*}
&&\left|E\left[\left(
\int_0^T f(r)\theta^{\e,k}(r)dr
\right)^{2m+1}\right]\right|\\
&\leq &\frac{1}{\e^{2m+1}} \int_{[0,T]^{2m+1}} \prod_{l=1}^{2m+1}|f(r_{l})|
\left|E\left[(-1)^{\sum_{i=1}^{2m+1} N(\frac{r_{i}}{\e^2})}\right]\right|
dr_1\ldots dr_{2m+1}\\
&=&\frac{2m+1}{\e^{2m+1}} \int_0^T |f(s)|ds\int_{[s,T]^{2m}}
\prod_{l=1}^{2m}|f(r_{l})|
\left|E\left[(-1)^{N(\frac{s}{\e^2})+\sum_{i=1}^{2m}
N(\frac{r_{i}}{\e^2})}\right]\right|
dr_1\ldots dr_{2m}\\
&\leq&(2m+1)\,\Delta^\e_{2m}(|f|,\ldots,|f|)\,\int_0^T
|f(s)|\frac{1}{\e}e^{-\frac{2s}{\e^2}}ds.
\end{eqnarray*}
Since for $s>\e$, we have that $\frac{1}{\e^2}{\rm
e}^{-\frac{2s}{\e^2}} \leq \frac12$, we get that
\begin{align}\label{inetwo}
&\int_0^T
|f(s)|\frac{1}{\e}e^{-\frac{2s}{\e^2}}ds
=\int_0^{\e} |f(s)|\frac{1}{\e}e^{-\frac{2s}{\e^2}}ds+ \e
\int_{\e}^T
|f(s)|\frac{1}{\e^2}e^{-\frac{2s}{\e^2}}ds\\
&\leq  \left( \int_0^{\e} |f(s)|^2 ds \right)^\frac12 \left(
\int_0^{\e} Q_\e(2s)ds
\right)^\frac12+ \frac{\e}{2} \int_{\e}^T |f(s)|ds
\leq \frac12 \left( \int_0^{\e} |f(s)|^2 ds \right)^\frac12 +
\frac{\e}{2} \|f\|_{L^2}, \nonumber
\end{align}
and (\ref{lm-cont-impair}) follows easily.

\end{proof}

The following lemma aims at measuring the distance between the laws of the stochastic integrals $\int_0^T f(r)\theta^{\e,k}(r)dr$ 
and $\int_0^T f(r)dW^k_r$, whenever $f$ is a given (deterministic) function: 
\begin{lemma}\label{lm-control2}
Let $f\in{\cal C}^\alpha([0,T])$ for a given $\alpha\in(0,1)$, $k\in\{1,2\}$
and $\e>0$.
For any $u\in\R$, we have:
\begin{eqnarray}
&&\left|E\big[e^{iu\int_0^T f(r)\theta^{\e,k}(r)dr}\big]
-E\big[e^{iu\int_0^T f(r)dW^k_r}\big]\right|\nonumber\\
&\leq& \left[\e^{2\alpha} c_\alpha\,\|f\|_\alpha\|f\|_{L^2} u^2 +
\phi_f(\e) \frac{u^2}{2} + \psi_f(\e) \frac{u^4}{8} +
\varphi_f(\e)\,\frac{|u|}{2}\right]{\rm
e}^{\frac{u^2\|f\|^2_{L^2}}{2}}, \label{lm-cont2}
\end{eqnarray}
with $c_\alpha=\int_0^\infty x^{\alpha}e^{-2x} dx$ and
$$
\phi_f(\e)=\int_0^T f^2(x)e^{-\frac{2x}{\e^2}}dx, \quad
\psi_f(\e)=\int_0^T dx \int_0^x  dy f^2(x) f^2(y){\rm
e}^{-\frac{2(x-y)}{\e^2}}.
$$
\end{lemma}
\begin{proof}
The proof is divided into two steps.\\
1. {\it First step: control of the imaginary part}. We can write, thanks to (\ref{lm-cont-impair}):
\begin{eqnarray*}
&&\left|{\rm Im}\left(E\big[e^{iu\int_0^T f(r)\theta^{\e,k}(r)dr}\big]
-E\big[e^{iu\int_0^T f(r)dW^k_r}\big]\right)\right|
=\left|{\rm Im}\,E\big[e^{iu\int_0^T f(r)\theta^{\e,k}(r)dr}\big]\right|\\
&\leq&\sum_{m=0}^\infty \frac{|u|^{2m+1}}{(2m+1)!}
\left|E\left[\left( \int_0^T f(r)\theta^{\e,k}(r)dr \right)^{2m+1}
\right]\right|\leq \varphi_f(\e)\,\frac{|u|}{2}\,{\rm
e}^{\frac{u^2\|f\|^2_{L^2}}{2}}.
\end{eqnarray*}
2. {\it Second step: control of the real part}. This step is more technical, and we will mainly get a bound on the quantity $L_{m,\e}$ defined by:
$$
L_{m,\e}=\left|\frac{1}{(2m)!}\,\Delta^\e_{2m}(f,\ldots,f)-\frac{1}{2^m}\int_0^T f^2(s_1)ds_1\ldots\int_0^{s_{m-1}}f^2(s_m)ds_m
\right|.
$$
In order to express this quantity in a suitable way for estimations, notice that
$\int_0^\infty e^{-2s}ds=\frac{1}{2}$. We  can thus insert this term artificially in the multiple integrals involved in the computations of $E[e^{iu\int_0^T f(r)dW^k_r}]$. This gives:
\begin{multline*}
L_{m,\e}=
\Bigg|\frac{1}{(2m)!}\,\Delta^\e_{2m}(f,\ldots,f)  \\
-\int_0^T f^2(r_1)dr_1\int_0^{\infty}e^{-2r_2}dr_2
\ldots\int_0^{r_{2m-3}}\!\!\!\!\!\!\!\!\!\!\!f^2(r_{2m-1})dr_{2m-1}\int_0^{\infty}e^{-2r_{2m}}dr_{2m}
\Bigg|.
\end{multline*}
By a telescoping sum argument, we can now write $L_{m,\e}$ as a sum of $m$ terms, whose prototype is given by 
$M_{m,\e}=M_{m,\e}^{1}+M_{m,\e}^{2}-M_{m,\e}^{3}$, with
\begin{multline*}
M_{m,\e}^{1}=
\int_0^T \!\! |f(r_1)|dr_1\int_0^{r_1} \!\!|f(r_2)| Q_\e(r_1-r_2) \, dr_2
\ldots\\
 \ldots\int_0^{r_{2m-2}}
\!\!\!\!\!\!|f(r_{2m-1})|dr_{2m-1}\int_0^{r_{2m-1}}\!\!\!\!\!\!
|f(r_{2m})-f(r_{2m-1})|Q_\e(r_{2m-1}-r_{2m}) \, dr_{2m},
\end{multline*}
where $M_{m,\e}^{2}$ is defined by
\begin{multline*}
M_{m,\e}^{2}=\int_0^T \!\!f(r_1)dr_1\int_0^{r_1} \!\!f(r_2)Q_\e(r_1-r_2)dr_2 \ldots\\
\ldots \int_0^{r_{2m-2}}
\!\!\!\!\!\!f^2(r_{2m-1}) dr_{2m-1} \int_{r_{2m-1}}^\infty
Q_\e(r_{2m})\, dr_{2m},
\end{multline*}
and where
\begin{equation*}
M_{m,\e}^{3}=
\int_0^T \!\!f(r_1)dr_1\int_0^{r_1} \!\!f(r_2)Q_\e(r_1-r_2)\,dr_2
\ldots\int_{r_{2m-2}}^{r_{2m-3}}
\!\!\!\!\!\!f^2(r_{2m-1})dr_{2m-1}\int_0^{\infty} \! e^{-2
r_{2m}}
dr_{2m}.
\end{equation*}
We will now bound those three terms separately: invoking first (\ref{ineone}), we get
\begin{eqnarray*}
M_{m,\e}^{1}
&\leq & \frac{1}{(m-1)!} \left( \int_0^T \!\!
|f(r_1)|dr_1\int_0^{r_1} \!\!|f(r_2)| Q_\e(r_1-r_2) \, dr_2 \right)^{m-1}\\
&&\quad\quad\quad \times \int_0^{T}
\!\!\!\!\!\!|f(r_{2m-1})|dr_{2m-1} \Vert f \Vert_\alpha
\int_0^{r_{2m-1}}\!\!\!\!\!\! |r_{2m}-r_{2m-1}|^\alpha  Q_\e(r_{2m-1}-r_{2m})\, dr_{2m}\\
&\leq & \frac{1}{(m-1)! 2^{m-1}} \Vert f \Vert_{L^2}^{2m-1} \Vert f
\Vert_\alpha c_\alpha \e^{2 \alpha}.
\end{eqnarray*}
On the other hand, (\ref{ineone}) and (\ref{inetwo}) also yield:
\begin{multline*}
M_{m,\e}^{2} \le \frac{1}{(m-1)!} \left( \int_0^T \!\! |f(r_1)|dr_1\int_0^{r_1}
\!\!|f(r_2)| Q_\e(r_1-r_2) \, dr_2
\right)^{m-1}  \\
\times\frac{1}{2} \int_0^{T} \!\!\!\!\!\! f^2(r_{2m-1})
Q_\e(r_{2m-1})\,    dr_{2m-1}
\leq  \frac{1}{(m-1)! 2^{m}} \Vert f \Vert_{L^2}^{2m-2}  \phi_{f}
(\e).
\end{multline*}
Finally, $M_{m,\e}^{3}$ can be bounded in the following way: observe that
\begin{multline*}
M_{m,\e}^{3}\le
\frac{1}{(m-2)!} \left( \int_0^T \!\! |f(r_1)|dr_1\int_0^{r_1}
\!\!|f(r_2)|Q_\e(r_1-r_2) \, dr_2 \right)^{m-2}\\
\times \frac12 \int_0^{T} \!\!\!\!\!\!|f(r_{2m-3})|dr_{2m-3}
\int_0^{r_{2m-3}}\!\!\!\!\!\! |f(r_{2m-2})| Q_\e(r_{2m-3}-r_{2m-2})\, dr_{2m-2}
\int_{r_{2m-2}}^{r_{2m-3}}\!\!\!\!\!\! f^2(r_{2m-1})dr_{2m-1},
\end{multline*}
which can also be written as
$$
M_{m,\e}^{3}\le \frac{1}{(m-2)! 2^{m-1}} \Vert f \Vert_{L^2}^{2m-4} M_{m,\e}^{3,1},
$$
with
\begin{eqnarray*}
M_{m,\e}^{3,1}&=&\int_0^T dr_{2m-1}f^2(r_{2m-1})
\int_{0}^{r_{2m-1}} \!\!\!\!\!\!dr_{2m-2}
\int_{r_{2m-1}}^T\!\!\!\!\!\! dr_{2m-3} |f(r_{2m-2})|
|f(r_{2m-3})| \\
&&\hskip7cm\times\frac{e^{-\frac{2(r_{2m-3}-r_{2m-1})}
{\e^2}}}{\e} \frac{e^{-\frac{2(r_{2m-1}-r_{2m-2})}
{\e^2}}}{\e}.
\end{eqnarray*}
It is now readily checked that
\begin{align*}
&M_{m,\e}^{3,1}\leq 
\int_{0}^{T}\!\!\!\!\!\! f^2(r_{2m-1})dr_{2m-1}
\int_{0}^{r_{2m-1}} \!\!\!\!\!\!dr_{2m-2}
\int_{r_{2m-1}}^T\!\!\!\!\!\! dr_{2m-3}
\\
 &\hspace{3cm}\times  \frac{ \left( f^2 (r_{2m-2})+
f^2(r_{2m-3}) \right)}{2} \frac{{\rm
e}^{-\frac{2(r_{2m-3}-r_{2m-1})} {\e^2}}}{\e} \frac{{\rm
e}^{-\frac{2(r_{2m-1}-r_{2m-2})} {\e^2}}}{\e}\\
&\leq  \frac{1}{2} 
\int_{0}^{T}\!\!\!\!\!\! f^2(r_{2m-1})dr_{2m-1}
\\
& \hspace{0.5cm} \times  \left( \int_{0}^{r_{2m-1}} \!\!\!\!\!\!dr_{2m-2}
 f^2
(r_{2m-2})  e^{-\frac{2(r_{2m-1}-r_{2m-2})} {\e^2}} +
\int_{r_{2m-1}}^T\!\!\!\!\!\! dr_{2m-3} f^2(r_{2m-3}) {\rm
e}^{-\frac{2(r_{2m-3}-r_{2m-1})}
{\e^2}} \right)\\
&\leq  \frac{1}{2}   \psi_f
(\e).
\end{align*}
Consequently,
$$
M_{m,\e}^{3}\le \frac{1}{(m-2)! 2^{m}} \Vert f \Vert_{L^2}^{2m-4} \psi_f.
$$

Our proof is now easily finished: plug our estimates on $M_{m,\e}^{1},M_{m,\e}^{2}$ and $M_{m,\e}^{3}$ into the definition of $M_{m,\e}$, and then in the definition $L_{m,\e}$. This yields
\begin{eqnarray*}
&&\left|{\rm Re}\left(E\big[e^{iu\int_0^T f(r)\theta^{\e,k}(r)dr}\big]
-E\big[e^{iu\int_0^T f(r)dW^k_r}\big]\right)\right|\\
&\leq&\sum_{m=1}^\infty \frac{u^{2m}}{(2m)!}\left|
\Delta_{2m}^\e (f,\ldots,f) - \frac{(2m)!}{2^mm!} \left(\int_0^T f^2(r)dr\right)^m
\right|\\
&\leq&\left[\e^{2\alpha} c_\alpha\,\|f\|_\alpha\|f\|_{L^2} u^2 +
\phi_f(\e) \frac{u^2}{2} + \psi_f(\e) \frac{u^4}{8}\right]{\rm
e}^{\frac{u^2\|f\|^2_{L^2}}{2}},
\end{eqnarray*}
which is our claim.

\end{proof}

The following lemma gives an alternative form for ${\bf X}^{\bf 2,\e}$ and ${\bf B}^{\bf 2}$:
\begin{lemma}\label{otherform}
Fix $i,j\in\{1,\ldots,d\}$, and $t>s\geq 0$. For all $\e>0$, we have
\begin{align}
&{\bf X}^{\bf 2,\e}_{st}(i,j) \notag\\
&=\int_0^t  \big(X_u^{j,\e}-X_s^{j,\e}\big) (t+\e-u)^{H-\frac12} \theta^{\e,i}(u)du
-\int_0^s  \big(X_u^{j,\e}-X_s^{j,\e}\big) (s+\e-u)^{H-\frac12} \theta^{\e,i}(u)du \notag\\
&-\al_H \int_0^t dv \theta^{\e,i}(v) \int_{s\vee v}^t du
(X_u^{j,\e}-X_v^{j,\e}) (u+\e-v)^{H-\frac32},
\label{lm2}
\end{align}
where we have set $\al_H=1/2-H$. In the limit $\e\to 0$, we also have
\begin{align}
&{\bf B}^{\bf 2}_{st}(i,j) \notag\\
&=\int_s^t \big(B^j_u-B^j_s\big)(t-u)^{H-\frac12}dW^i_u
-\int_0^s \big(B^j_u-B^j_s)\big[(t-u)^{H-\frac12}-(s-u)^{H-\frac12}\big]dW^i_u\notag\\
&-\al_H\int_0^t dW^i_v\int_{v\vee s}^t du\big(B^j_u-B^j_v\big)(u-v)^{H-\frac32}.
\label{lm12}
\end{align}
\end{lemma}
\begin{proof}
For any $\e>0$, the process $X^{\e,i}$ is differentiable, and according to (\ref{eq:def2-x-ep}), we have
$$
\dot X^{\e,i}(r)=\e^{H-1/2} \theta^{\e,i}(u) 
-\al_H \int_0^u (u+\e-v)^{H-3/2} \theta^{\e,i}(v) \, dv.
$$
Recall also that we have set $\der X^{j,\e}_{st}=X^{j,\e}_{t}-X^{j,\e}_{s}$ for any $s,t\in\ott$. This allows to write:
\begin{multline}\label{eq:54}
{\bf X}^{\bf 2,\e}_{st}(i,j) =\int_s^t  \der X^{j,\e}_{su}\, dX_u^{i,\e}  
=\e^{H-\frac12}\int_s^t  \der X^{j,\e}_{su}\, \theta^{\e,i}(u)du  \\
-\al_H \int_s^t  du \, \der X^{j,\e}_{su}
\int_0^u dv (u+\e-v)^{H-\frac32}\theta^{\e,i}(v).
\end{multline}
Moreover, an elementary application of Fubini's theorem yields:
\begin{align*}
&\int_s^t  du \, \der X^{j,\e}_{su}
\int_0^u dv (u+\e-v)^{H-\frac32}\theta^{\e,i}(v)  
=
\int_{0}^t dv\theta^{\e,i}(v) 
\int_{s\vee v}^t du \, \der X^{j,\e}_{su} \, (u+\e-v)^{H-\frac32}  \\
&=\int_0^t dv \theta^{\e,i}(v) \, \delta X^{j,\e}_{sv} 
\int_{s\vee v}^t dr  (r+\e-u)^{H-\frac32} 
+\int_0^t dv \theta^{\e,i}(v) \int_{s\vee v}^t dr \,
\delta X^{j,\e}_{vr}  
(r+\e-v)^{H-\frac32}
\end{align*}
Integrating the kernel $(r+\e-u)^{H-\frac32}$, and plugging the last identity into (\ref{eq:54}), we obtain the desired relation (\ref{lm2}).

\smallskip

To get  formula (\ref{lm12}) for ${\bf B}^{\bf 2}_{st}(i,j)$, it suffices to observe
that $${\bf B}^{\bf 2}_{st}(i,j)=L^2-\lim_{\e\to 0}\int_s^t \big(B^j_u-B^j_s)dB^{i,\e}_u$$
with $B^{i,\e}_u=\int_0^u (u+\e-v)^{H-\frac12}dW^i_v$,
and then to
mimick the computations allowing us to write (\ref{lm2}) just above.
Details are left to the reader (see also the proof of \cite[Lemma 3]{AMN}).

\end{proof}

Finally, the following lemma gives an estimate for the variance of ${\bf B}^{\bf 2}_{st}(i,j)$
which is useful in the proof of Proposition \ref{prop:hyp-fbm}:
\begin{lemma}\label{lm-technical}
There exists a constant $c>0$, depending only on $H$, such
that $E\big|{\bf B}^{\bf 2}_{st}(i,j)\big|^2$ $\leq c|t-s|^{4H}$ for all $t>s\geq 0$
and $i,j\in\{1,\ldots,d\}$.
\end{lemma}
\begin{proof}
The case where $i=j$ is immediate by Lemma \ref{onenabesoin}, so we only concentrate on the case where
$i\neq j$.
Using formula (\ref{lm12}), we see that we have
to bound the three following terms:
\begin{eqnarray*}
A_1&:=&\int_s^t E\big|B^j_u-B^j_s\big|^2(t-u)^{2H-1}du\\
A_2&:=&\int_0^s  E\big|B^j_u-B^j_s\big|^2\left((t-u)^{H-\frac12}-(s-u)^{H-\frac12}\right)^2 du\\
A_3&:=&\int_0^t E\left|\int_{v\vee s}^t du\big(B^j_u-B^j_v\big)(u-v)^{H-\frac32}\right|^2 dv.
\end{eqnarray*}
Throughout the proof, $c$ will denote a generic constant (depending only on $H$, $T$) whose
value can change from one line to another.
Owing to the fact that $E\big|B^j_u-B^j_s\big|^2\leq c|u-s|^{2H}$, see Lemma \ref{onenabesoin}, we can write
$$
A_1\leq c\int_s^t (u-s)^{2H}(t-u)^{2H-1}du\leq c(t-s)^{2H}\int_s^t (t-u)^{2H-1} = c(t-s)^{4H}.
$$
We also get
\begin{eqnarray*}
A_2 &\leq& c\int_0^s  (s-u)^{2H}
\left((t-u)^{H-\frac12}-(s-u)^{H-\frac12}\right)^2 du\\
&=&c\int_0^s  u^{2H}
\left((t-s+u)^{H-\frac12}-u^{H-\frac12}\right)^2 du\\
&=&c(t-s)^{4H}\int_0^{\frac{s}{t-s}}  u^{2H}
\left((1+u)^{H-\frac12}-u^{H-\frac12}\right)^2 du\\
&\leq&c(t-s)^{4H}\int_0^{\infty}  u^{2H}
\left((1+u)^{H-\frac12}-u^{H-\frac12}\right)^2 du
= c(t-s)^{4H},
\end{eqnarray*}
the last integral being finite since $H<\frac12$.
Finally, we have
\begin{eqnarray*}
&&E\left|\int_{v\vee s}^t du\big(B^j_u-B^j_v\big)(u-v)^{H-\frac32}\right|^2 \\
&=&\int_{v\vee s}^t du \int_{v\vee s}^t dw E\big[(B^j_u-B^j_v)(B^j_w-B^j_v)\big] (u-v)^{H-\frac32}(w-v)^{H-\frac32}\\
&\leq& c\int_{v\vee s}^t du \int_{v\vee s}^t dw (u-v)^{2H-\frac32}(w-v)^{2H-\frac32}
=c\left(\int_{v\vee s}^t (u-v)^{2H-\frac32} du\right)^2
\end{eqnarray*}
so that
\begin{eqnarray*}
A_3 &\leq&c \int_0^s \big[ (t-v)^{2H-\frac12}-(s-v)^{2H-\frac12}  \big]^2 dv
+c\int_s^t \left(\int_v^t (u-v)^{2H-\frac32} du\right)^2 dv\\
&\leq& c (t-s)^{4H}\int_0^\infty \big[ (1+v)^{2H-\frac12}-v^{2H-\frac12}  \big]^2 dv + c\int_s^t (t-v)^{4H-1} dv\\
&=& c(t-s)^{4H}. 
\end{eqnarray*}
This finishes the proof of the lemma.

\end{proof}


\end{document}